%% file: inscribed_curvature_paper.tex
\documentclass[11pt,reqno]{amsart}
\usepackage[T1]{fontenc}
\usepackage[utf8]{inputenc}
\usepackage[english]{babel}
\usepackage{microtype}

\usepackage{lmodern}

\usepackage[a4paper,margin=1in]{geometry}

\usepackage{amsmath,amssymb,amsthm,mathtools}
\mathtoolsset{showonlyrefs} 
\usepackage{mathtools} 

\usepackage{graphicx}
\usepackage{tikz}
\usepackage[labelfont=bf]{caption}
\usepackage{subcaption}
\usepackage{enumitem}
\usetikzlibrary{backgrounds}

\usepackage{csquotes}
\usepackage{xcolor}
\usepackage[backend=biber,style=alphabetic,maxbibnames=9]{biblatex}
\usepackage[colorlinks=true,linkcolor=blue!60!black,citecolor=green!40!black,urlcolor=blue!60!black]{hyperref}
\usepackage[nameinlink,noabbrev,capitalise]{cleveref}

\numberwithin{equation}{section}
\theoremstyle{plain}
\newtheorem{theorem}{Theorem}[section]
\newtheorem{lemma}[theorem]{Lemma}
\newtheorem{proposition}[theorem]{Proposition}
\newtheorem{claim}[theorem]{Claim}
\newtheorem{corollary}[theorem]{Corollary}
\theoremstyle{definition}
\newtheorem{definition}[theorem]{Definition}
\theoremstyle{remark}
\newtheorem{remark}[theorem]{Remark}
\newtheorem{example}[theorem]{Example}

\crefname{claim}{Claim}{Claims}
\crefname{example}{Example}{Examples}
\AddToHook{env/lemma/begin}{\crefalias{theorem}{lemma}}
\AddToHook{env/proposition/begin}{\crefalias{theorem}{proposition}}
\AddToHook{env/claim/begin}{\crefalias{theorem}{claim}}
\AddToHook{env/corollary/begin}{\crefalias{theorem}{corollary}}
\AddToHook{env/definition/begin}{\crefalias{theorem}{definition}}
\AddToHook{env/remark/begin}{\crefalias{theorem}{remark}}
\AddToHook{env/example/begin}{\crefalias{theorem}{example}}

\newcommand{\RR}{\mathbb{R}}

\newcommand{\Hess}{\operatorname{Hess}}
\newcommand{\Rm}{\operatorname{Rm}}

\newcommand{\seq}{\subseteq}
\newcommand{\ov}{\overline}

\DeclareMathOperator{\Ric}{Ric}
\DeclareMathOperator{\inj}{inj}
\DeclareMathOperator{\Span}{Span}

\DeclareMathOperator{\Diam}{Diam}
\DeclareMathOperator{\Id}{Id}
\DeclareMathOperator{\osc}{osc}
\DeclarePairedDelimiter{\ip}{\langle}{\rangle}

\usepackage{footmisc}
\title[A Chord-Arc Estimate for Graphs]{A Chord-Arc Estimate for
Graphical Hypersurfaces}
\author[Jacob Reznikov]{Jacob Reznikov}
\address{Department of Mathematics, Massachusetts Institute of
Technology, Cambridge, MA, USA}
\email{jacob\_oc@mit.edu}
\date{November 10th, 2026}

\begin{document}
\begin{abstract}
  We prove an interior chord-arc estimate for graphical hypersurfaces without an a priori gradient bound, which is not restricted to balls on the scale of their curvature. As a main application we derive an in-radius lower bound depending on the curvature which then also implies Heinz-type estimates for minimal graphs and uniform area bounds for graphical hypersurfaces.
\end{abstract}


\maketitle

\input{sections/introduction}
\input{sections/curvature_scale_hypersurfaces}
\input{sections/chord_arc_bound}
\input{sections/inscribed_curvature}

\input{sections/main_theorem}
\input{sections/examples}
\input{sections/applications.tex}
\input{sections/app_differential_geometry}
\input{sections/app_analytic_estimates}
\input{sections/app_geometric_identities}

\printbibliography

\end{document}

%% file: sections/introduction.tex
\section{Introduction}

A chord-arc estimate bounds the intrinsic distance between two points of a
surface by their extrinsic distance.
Surfaces with small chord-arc constant were introduced by Semmes
\cite{semmesChordArcI1991,semmesChordArcII1991}, as a higher-dimensional
version of the planar chord-arc curves studied by David
\cite{davidCourbesCordeArc1982}. Semmes considered complete hypersurfaces
$\Gamma \subset \RR^{n+1}$ that look flat in every ball, in the sense that the
unit normal has small BMO norm and $\Gamma$ stays close to a hyperplane, and
he showed that $d_\Gamma(x,y)/|x-y|$ is then close to $1$ for all
$x,y \in \Gamma$.

Chord-arc bounds later became an important tool in minimal surface theory. Colding and Minicozzi proved one for
embedded minimal disks in $\RR^3$ \cite{coldingMinicozziCalabiYau2008}, and it
was the main geometric step in their proof of the Calabi--Yau conjectures for
embedded surfaces. Meeks and Tinaglia extended their bound to embedded disks
of constant mean curvature \cite{meeksTinagliaChordArc2018}.

Our estimate is of Semmes's kind, with ratio close to $1$, but it is for
graphical hypersurfaces in a product $N \times \RR$ and the hypotheses are
different. Our only assumptions are a bound on the second fundamental form, and that the graph has to be
defined over a thin neighborhood of a geodesic of $N$.

There is a trivial bound one gets from a bounded gradient, but one should expect a better bound in terms of curvature. For instance a highly sloped plane has a chord-arc ratio of exactly $1$ whereas the naive gradient estimate can be made arbitrarily large. Our estimate also differs from most estimates in the literature in that it is only local in the base: the two points in question can be chosen arbitrarily far away in the ambient space.

From now on, unless stated otherwise, $N$ denotes a
complete Riemannian manifold with bounded curvature, $\Omega \subset N$ an
open set, and $\Sigma \subset N \times \RR$ a graphical
hypersurface over $\Omega$.
For $x \in \Omega$, write $\hat{x}$ for its lift to $\Sigma$.
For $x$ in $N$ or $N\times\RR$ (after projecting to $N$ if needed), write
\mathtoolsset{showonlyrefs=false}
\begin{align}
	\rho(x) & := \min(d(x, \partial \Omega), \inj(N)), \label{eq-rho} \\
	\kappa  & := \|h_\Sigma\|_{C^0(\Sigma)}, \label{eq-kappa}         \\
	K       & := \|\Rm_N\|_{C^0}. \label{eq-ambient-curvature}
\end{align}
\mathtoolsset{showonlyrefs=true}

Fix points $x,y \in \Omega$, set $r = d(x,y)$, and assume $\gamma$ is an
interior minimizing geodesic of $N$ from $x$ to $y$.

\begin{theorem}[Chord-arc estimate]\label{lem:intro-chord-arc}
	Given a scale $\varepsilon \leq r$ such that the $\varepsilon$-neighborhood
	$B_\varepsilon(\gamma)$ of $\gamma$ lies in $\Omega$, define the
	scale-invariant quantity $\delta$ as
	\[
		\delta := r^2 \left(\kappa\varepsilon^{-1} + K\right).
	\]
	If $\delta$ is sufficiently small, then
	\[
		d_\Sigma(\hat{x},\hat{y})
		\leq (1+C\delta)\,d_{N\times\RR}(\hat{x},\hat{y}),
	\]
	where $C$ is dimensional.
\end{theorem}

The width $\varepsilon$ cannot be removed from $\delta$:
in \cref{ex-cylinder} we have $r = 1$ and the ratio of the two distances is
at least $1 + c\,\kappa/\varepsilon$ for a dimensional constant $c > 0$.
On the other hand, $d_{N\times\RR}(\hat{x},\hat{y})$ is not bounded in terms of
$r$. 

The main application of \cref{lem:intro-chord-arc} is a lower bound on the
in-radius.
We write $r_*(p)$ for the radius of the largest one-sided tangent ball at $p$,
the in-radius of $\Sigma$ at $p$ (\cref{def-in-radius}). It depends on the
choice of orientation of $\Sigma$.
For embedded hypersurfaces a curvature bound does not bound $r_*$ from below:
the hourglass surfaces $\Sigma_R$ of \cref{ex-hourglass} have
$|A_{\Sigma_R}| \leq C/R$ but $r_* \leq 1$ at $(0,1)$, and converge to two
parallel planes as $R \to \infty$. For graphs it does, up to ambient and
boundary terms.

\begin{theorem}[In-radius bound]\label{thm:main}
	There exists a dimensional constant $C > 0$ such that for either orientation
	of $\Sigma$ and all interior points $p \in \Sigma$,
	\[
		r_*(p)^{-1} \leq C \bigl(\kappa + K^{1/2}\bigr)
		+ 4 \rho(p)^{-1}.
	\]
\end{theorem}

In the proof of \cref{thm:main}, the maximal tangent ball at $p$ touches a
graph at a second point $q$ from the same side, so the normal turns by about
$d(p,q)/r_*$ from $p$ to $q$, while \cref{lem:intro-chord-arc} gives a path
from $p$ to $q$ of length comparable to $d(p,q)$, along which it turns by at
most $\kappa$ times the length. These are incompatible when $\kappa r_*$ is
small.

\begin{remark}
	In a Riemannian ambient space one can also measure non-collapsing by an
	inscribed curvature form $h_\Sigma^*$, built from the second fundamental
	forms of the admissible tangent balls at $p$, in the spirit of
	\cite{andrewsHanLiWei2015}. Since $h_\Sigma^* \geq h_\Sigma$, and by Hessian
	comparison, $-\kappa g \leq h_\Sigma^* \leq (r_*^{-1} + K^{1/2})g$, so \cref{thm:main}
	bounds $\|h_\Sigma^*\|$ as well, with a larger constant $C$.
\end{remark}

\begin{remark}
	There is also a different approach to inscribed-curvature estimates which
	does not use such a chord-arc input; it will be treated in an upcoming work.
\end{remark}

\subsection*{Examples and applications}
Section~\ref{sec:examples} shows by examples that graphicality and
the scale in \cref{lem:intro-chord-arc} are needed.
Section~\ref{sec:apps} gives a Heinz-type estimate for the
in-radius of minimal graphs, Corollary~\ref{cor-heinz-inscribed-curvature},
and uniform area estimates for graphical hypersurfaces in curved ambient
spaces, Corollaries~\ref{cor-local-area-estimate} and \ref{cor-uniform-area-estimate}.

\subsection*{Related work}
The in-radius goes back to Blaschke’s Rolling Theorem \cite[\S24]{blaschkeKreisKugel1916}
and was more recently used by Andrews \cite{andrewsNoncollapsingMeanconvexMean2011}
in mean curvature flow, and by Andrews--Langford--McCoy
\cite{andrewsNoncollapsingFullyNonlinear2013} and Andrews--Langford
\cite{andrewsTwosidedNoncollapsingCurvature2014} in more general curvature flows.
Lower bounds for the normal injectivity radius, or reach
\cite{federerCurvatureMeasures1959}, need a global hypothesis besides a bound
on the second fundamental form. Boldt, G\"uneysu and Pigola
\cite{boldtGuneysuPigola2025} use a topological condition on small balls.
For graphs, \cref{thm:main} bounds both in-radii from below without such a
condition or a gradient bound.
Uniform regular neighborhoods were constructed for properly embedded
minimal surfaces with bounded curvature in $\RR^3$
\cite{meeksRosenbergMaxPrinciples2008} and for strongly Alexandrov embedded
hypersurfaces of constant mean curvature $H > 0$
\cite{meeksTinagliaRegularNeighborhoods2011}.

From this point onward, we allow $c,C$ to vary from line to line.
They denote positive dimensional constants, with $C$ chosen sufficiently
large and $c$ sufficiently small. Constants with additional dependencies
are written explicitly.

\subsection*{Acknowledgments}
I would like to thank my advisor, Prof. W. Minicozzi, for the project that inspired this paper, and for his support and helpful discussions during its writing.
During this project, I was partially supported by NSF DMS Grant 2304684 and by the Natural Sciences and Engineering Research Council of Canada (NSERC) [grant number CGRSD-609347-2026].

\subsection*{AI Disclosure}
Apart from the proof of \cref{prop-connection-expansion}, no AI was used for the mathematics in this paper. AI tools were used for editing and for searching for related work.

%% file: sections/curvature_scale_hypersurfaces.tex
\section{Hypersurfaces at the Curvature Scale}\label{sec:curvature-scale}
The goal of this section is to establish quantitative properties for
hypersurfaces with a fixed curvature bound. These estimates are scale-covariant:
if $\|A_\Sigma\| < \kappa$, then the relevant geometric scale is
$\kappa^{-1}$, and the estimates are controlled by quantities such as
$\kappa r$, $\kappa t^2$, and $\kappa R$. They are used to prove the chord-arc
bound, where controlled curvature at the relevant scale must be converted into
control of geodesics, graphical representations, and nearly-parallel level sets.
For the rest of the section we fix a surface $\Sigma$ embedded
into $\RR^{n+1}$ such that $\|A_\Sigma\| < \kappa$ everywhere. We allow
boundary, so for a point $x$ we let $d(x, \partial \Sigma)$ denote the
distance to the boundary, or $\infty$ if the boundary is empty.
For any ambient domain $X \subseteq \RR^{n+1}$ we denote by
$X_\varepsilon$ the set $\{x \in X : d(x, \partial X) \geq
	\varepsilon \}$.
Balls written without a center, such as $B_s$ or $B_s^P$, are centered at the
origin.

\begin{proposition}\label{prop-curvature-scale-basic}
	For a hypersurface $\Sigma$ satisfying the curvature bound above, the
	following holds:
	\begin{enumerate}[label=(\roman*)]
		\item The normal vector map $\Sigma \to S^n$ is $\kappa$-Lipschitz.
		\item Geodesics $\gamma : [0,T] \to \Sigma$ satisfy
		      $\|\partial_t^2 \gamma\| < \kappa$ as curves in $\RR^{n+1}$.
		      Furthermore, for $0 \leq t \leq T$,
		      \begin{align*}
			      \|\gamma(t) - \gamma(0) - t\gamma'(0)\| & \leq \frac{\kappa t^2}{2},
			      \\
			      \|\gamma(t) - \gamma(0)\|               & \geq t\left(1-\frac{\kappa^2 t^2}{6}\right).
		      \end{align*}
		\item For affine subspace $S$ that intersects $\Sigma$ transversally we
		      have that $\|A_{\Sigma \cap S}\| < \kappa \|\pi_S
			      \nu_\Sigma\|^{-1}$.
	\end{enumerate}
\end{proposition}
\begin{proof}
	$(i)$ is trivial, since $\nabla \nu$ is exactly the second
	fundamental form. $(ii)$ follows quickly, since $\gamma$ is an
	intrinsic geodesic we have
	\[
		\nabla^{\RR^{n+1}}_{\partial_t \gamma} (\partial_t \gamma)
		= \nabla^{\Sigma}_{\partial_t \gamma} (\partial_t \gamma) +
		(\nabla^{\RR^{n+1}}_{\partial_t \gamma} - \nabla^{\Sigma}_{\partial_t
			\gamma})(\partial_t \gamma)
		= 0 - h(\partial_t \gamma, \partial_t \gamma) \nu,
	\]
	and this is bounded by $\kappa$ since $\partial_t \gamma$ is a unit vector.
	Integrating twice gives the Taylor remainder estimate. Also,
	$\|\gamma'(s)-\gamma'(0)\| \leq \kappa s$, so
	\[
		\ip{\gamma'(s),\gamma'(0)}
		= 1-\frac{\|\gamma'(s)-\gamma'(0)\|^2}{2}
		\geq 1-\frac{\kappa^2 s^2}{2}.
	\]
	Integrating from $0$ to $t$ and using
	$\|\gamma(t)-\gamma(0)\| \geq \ip{\gamma(t)-\gamma(0),\gamma'(0)}$
	gives the chord bound.
	For $(iii)$, first note that $\nu_{\Sigma \cap S} = \frac{\pi_S
			\nu_\Sigma}{\|\pi_S \nu_\Sigma\|}$ we pick any orthonormal frame $e_i$ for
	$\Sigma \cap S$ and then compute
	\[
		h_{\Sigma \cap S} (e_i,e_j)
		= \ip*{\nabla_{e_i} \nu_{\Sigma \cap S},e_j}
		= \ip*{\nabla_{e_i} \frac{\pi_S \nu_\Sigma}{\|\pi_S \nu_\Sigma\|},e_j}
		= \|\pi_S \nu_\Sigma\|^{-1} \ip{\pi_S \nabla_{e_i} \nu_\Sigma, e_j},
	\]
	and so since $e_i,e_j$ are unit vectors, the second term is bounded above
	by $\kappa$.
\end{proof}

\begin{corollary}\label{cor-gradient-oscilation}
	Let $\Sigma$ satisfy $\|A_\Sigma\| < \kappa$ and be graphical over a domain
	$B$. Then we have
	\[
		\left\|\nabla^{N} \frac{1}{1 + \|\nabla^N f\|_N^2}\right\| \leq 2\kappa,
	\]
	where $f$ is the function whose graph is $\Sigma$.
\end{corollary}
\begin{proof}
	Choose coordinates with $e_1$ vertical and $\ip{e_1,\nu}>0$.
	Let $w$ be an arbitrary unit vector at a point $p \in B$, and
	$\hat{e}_j$ be an orthonormal frame for $\Sigma$. Then we have
	\[
		\nabla_{w} \frac{1}{\sqrt{1+\|\nabla^N f\|_N^2}}
		= \nabla_{w} \ip{e_1, \nu}
		= \ip{e_1, \nabla_{w} \nu}
		= \ip{e_1, A(w)},
	\]
	where $A$ is the shape operator of $\Sigma$. The norm of this operator is
	$\kappa$ in the metric induced by $\Sigma$ by assumption, so we have
	$\|A(w)\| \leq \kappa \|w\|_\Sigma$. We know that in the coordinates of the
	base space, the metric has eigenvalues $1$ and $1 + \|\nabla^N f\|_N^2$,
	so we have $\|w\|_\Sigma \leq \sqrt{1 + \|\nabla^N f\|_N^2}$. Combining these
	we get
	\[
		\left\|\nabla_{w} \frac{1}{\sqrt{1+\|\nabla^N f\|_N^2}}\right\|
		\leq \kappa \sqrt{1 + \|\nabla^N f\|_N^2}.
	\]
	Now applying the power rule on the LHS we get the desired result.
\end{proof}

Especially useful is the second statement, which can be used to prove the
following.
\begin{proposition}\label{prop-geodesic-boundary-seperation}
	Assume that $\Sigma$ is properly embedded inside a convex domain
	$D$ with $\partial \Sigma \subseteq \partial D$. Then
	geodesics $\gamma : [0,T] \to \Sigma$ which are interior for $t > 0$,
	with $\gamma(0) \in \partial \Sigma$ and
	$\gamma'(0) \in T_{\gamma(0)} \partial \Sigma$, satisfy
	$d(\gamma(t), \partial D) \leq \kappa t^2/2 < C\kappa t^2$ for $t > 0$, thus
	$\gamma(t) \not \in D_{C\kappa t^2}$.
\end{proposition}
\begin{proof}
	Let $H$ be a supporting half space for $D$ at $\gamma(0)$. Since
	$\gamma'(0)$ is tangent to $\partial\Sigma \subseteq \partial D$, the line
	$\gamma(0)+t\gamma'(0)$ lies in $\partial H$. Thus
	\cref{prop-curvature-scale-basic} part $(ii)$ gives
	\[
		d(\gamma(t),\partial D) \leq d(\gamma(t),\partial H)
		\leq \|\gamma(t)-\gamma(0)-t\gamma'(0)\| \leq \frac{\kappa t^2}{2}.
	\]
\end{proof}

\begin{corollary}\label{cor-geodesically-connected}
	In the same setting as in \cref{prop-geodesic-boundary-seperation},
	let $p,q$ be points in $\Sigma \cap D_{C\kappa r^2}$ with $d^\Sigma
		(p, q) < r$. Then the length minimizing path between $p$ and $q$ is
	an interior geodesic.
\end{corollary}
\begin{proof}
	Let $\gamma$ be the length minimizing path between $p$ and $q$.
	Assume that $\gamma(t_0)$ is a point in $\partial \Sigma$ for some
	$0 < t_0 < r$, by closedness of the boundary we may assume that
	$\gamma(t)$ is in the interior for all $t < t_0$. By
	\cite[Theorem~1(A)]{alexanderGeodesicsRiemannianManifoldswithBoundary1981} we
	know that $\gamma$ is $C^1$ and so $\gamma'(t_0)$ is tangent to
	$\partial \Sigma$. Furthermore we know that $\gamma |_{[0,t_0]}$ is
	an interior geodesic. Hence applying
	\cref{prop-geodesic-boundary-seperation} to $\gamma |_{[0,t_0]}$ but
	reversed, we get
	$d(\gamma(0),\partial D)< C\kappa t_0^2\leq C\kappa r^2$. This
	contradicts our initial assumption, and hence no point of $\gamma$
	is contained in the boundary. It is thus an interior geodesic.
\end{proof}

It is known that all embedded submanifolds are locally graphs over their
tangent plane, but under a curvature bound at the scale $\kappa^{-1}$ we can
turn this qualitative result into a quantitative result.
\begin{proposition}\label{prop-qgraph}
	Let $\Sigma$ be a hypersurface with boundary, properly embedded in an open
	set $U \seq \RR^{n+1}$, which satisfies $\|A_\Sigma\| < \kappa$ everywhere. Let
	$p$ be a point on $\Sigma$ with tangent plane $P$. There are constants
	$c,C>0$ such that $B_{C r}^\Sigma (p)$ contains a graph of some function
	$u$ over $B_r^P(p)$ whenever
	$r \leq c \min\{\kappa^{-1}, d(p,\partial \Sigma), d(p,\partial U)\}$.

	Furthermore, the graph satisfies $\|\nabla u (x)\| \leq C \kappa \|x - p\|$ as
	well as $|u(x)| \leq C \kappa \|x - p\|^2$.
\end{proposition}
\begin{proof}
	We will assume WLOG that $p = 0$ and $P = \{0\} \times \RR^n$, so
	that $\nu(0) = e_1$. Then for $q \in B_{C r}^\Sigma (p)$, we have
	$\|\nu(q) - e_1\| < C \kappa r \leq c$, so by
	\cref{prop-simple-geometric} part $(i)$ we have $\ip{e_1, \nu} \geq c$.
	Now take an orthonormal frame $\hat{e}_i$ for $\Sigma$, and
	rotate it such that $\hat{e}_2,\dots,\hat{e}_n$ are all orthogonal
	to $e_1$. Denote by $\pi_P$ the projection onto $P$, then clearly
	$\pi_P \hat{e}_i = \hat{e}_i$ for $i \geq 2$ and
	\[
		\|\pi_P \hat{e}_1\| = \sqrt{1 - \ip{e_1,\hat{e}_1}^2} = |\ip{e_1,
			\nu(q)}| \geq c.
	\]
	Now applying \cref{prop-qinverse}, and increasing $C$ if necessary, we get
	that $\pi_P (B_{C r}^\Sigma)$ contains $B_{r}^P$.

	To show it is injective over this ball, assume that distinct points $x,y \in
		B_{C r}^\Sigma(p)$ have $\pi_P (x) = \pi_P (y)$. By properness there is a minimizing path $\gamma$ from $x$ to $y$ of length $\ell < 2Cr$.
	Every point of $\gamma$ is within distance $3Cr$ of $p$, while
	$d(p,\partial\Sigma)>3Cr$ and $d(p,\partial U)>3Cr$ after shrinking $c$.
	Thus $\gamma$ is an interior geodesic. Since $\|\nu(x)-e_1\|<c$, we have
	$\|\pi_P\gamma'(0)\|\geq 1/2$, and \cref{prop-curvature-scale-basic} part $(ii)$
	gives
	\[
		0 = \|\pi_P(y-x)\|
		\geq \ell\|\pi_P\gamma'(0)\|-\frac{\kappa\ell^2}{2}
		\geq \frac{\ell}{2}-\frac{\kappa\ell^2}{2}>0,
	\]
	a contradiction for $c$ sufficiently small. Hence the map is injective.

	Thus $\Sigma$ contains the graph of a function $u$ over $B_{r}^P$. Since
	for graphs we have $\nu_u = (1+\|\nabla u\|^2)^{-1/2} (1, -\nabla u)$, we
	get that $\ip{\nu_u, e_1} = (1+\|\nabla u\|^2)^{-1/2}$ and $\pi_P \nu_u = -
		(1+\|\nabla u\|^2)^{-1/2} \nabla u$. We combine these two facts to get
	\[
		\|\nabla u(x)\| = \frac{\|\pi_{P} \nu_u\|}{\ip{\nu_u, e_1}} \leq
		C \|\nu_u - e_1\| < C r \kappa.
	\]
	Now applying the same argument but with $r$ replaced with $\|x\|$ we get
	$\|\nabla u(x)\| < C \|x\| \kappa$ and hence we get our estimate on the
	derivative. Finally integrating this estimate we get $|u(x)| \leq C \kappa
		\|x\|^2$.
\end{proof}

We end this section with a key property about hypersurfaces with the same
curvature bound which are close but disjoint on the scale of their curvature.
Namely, they must be nearly-parallel.
\begin{proposition}\label{prop-disjoint-curvature-scale}
	Fix scale $R \leq c\kappa^{-1}$. Let $\Sigma_1,\Sigma_2$ be two
	connected hypersurfaces satisfying $\|A_{\Sigma_i}\| < \kappa$ and properly
	embedded inside $B_{C R}$. If
	the surfaces are disjoint, and close in the sense that both
	intersect $B_{r}$ for some positive $r < cR$, while
	$\partial \Sigma_i \cap B_{C R} = \emptyset$, then there are choices of
	unit normal fields on $\Sigma_1$ and $\Sigma_2$ such that
	$\|\nu_{\Sigma_1}(p) - \nu_{\Sigma_2}(q)\| < C r R^{-1} + C \kappa
		R$ for all $p \in \Sigma_1, q \in \Sigma_2$.
\end{proposition}
\begin{proof}
	For $c$ small, $f=\|x\|^2$ satisfies
	$\Hess_{\Sigma_i}f\geq 2(1-\kappa CR)g_{\Sigma_i}\geq g_{\Sigma_i}$.
	Properness and connectedness place any two points inside a compact regular
	sublevel component of $f$. Its strictly convex boundary keeps a minimizing
	path $\gamma$ interior. If $\gamma$ has length $\ell$, then
	$0\leq f(\gamma(\ell/2))\leq (CR)^2-\ell^2/8$, giving
	$\Diam_{\Sigma_i}(\Sigma_i)<3CR$. By
	\cref{prop-curvature-scale-basic} part $(i)$ we get
	\[
		\osc_{\Sigma_i} \nu_{\Sigma_i} < C \kappa R.
	\]

	Pick any $p_0 \in \Sigma_1 \cap B_r$, and pick $q_0 \in \Sigma_2$ with minimal distance to
	$p_0$. Since $\Sigma_2$ intersects $B_r$, we get that
	$a := \|p_0 - q_0\| < C r$ and $q_0 \in B_{C r}$. The diameter bound gives
	\[
		d^{\Sigma_1}(p,p_0) < 3 C R, \qquad
		d^{\Sigma_2}(q,q_0) < 3 C R,
	\]
	for all $p \in \Sigma_1$ and $q \in \Sigma_2$. In particular,
	\[
		\|\nu_{\Sigma_1}(p) - \nu_{\Sigma_1}(p_0)\| < C \kappa R, \qquad
		\|\nu_{\Sigma_2}(q) - \nu_{\Sigma_2}(q_0)\| < C \kappa R.
	\]

	We now estimate the normals at the single pair $p_0,q_0$. Since
	$p_0 \in B_r$, $q_0 \in B_{C r}$ and $r<cR$, the assumption
	$\partial\Sigma_i \cap B_{C R} = \emptyset$ gives
	$d(p_0,\partial\Sigma_1)\geq cR$ and $d(q_0,\partial\Sigma_2)\geq cR$.
	Thus $q_0$ is an interior point of $\Sigma_2$ at which $\|p_0-\cdot\|^2$
	attains its minimum on $\Sigma_2$, so $p_0-q_0$ is normal to $\Sigma_2$ at
	$q_0$. As the surfaces are disjoint, $a>0$, and we orient $\Sigma_2$ so
	that $\nu_{\Sigma_2}(q_0)=\nu$, where $\nu = (p_0-q_0)/a$. Let
	$P = q_0 + \nu^\perp$ be the tangent plane of $\Sigma_2$ at $q_0$, and
	$\pi_P$ the orthogonal projection onto $P$.

	By \cref{prop-qgraph}, applied at $q_0$ with $U = B_{CR}$ at scale $cR$,
	$\Sigma_2$ contains the graph $\{x+u(x)\nu : x\in B_{cR}^{P}(q_0)\}$ of a
	function $u$ over $B_{cR}^{P}(q_0)$ satisfying
	$|u (x)| \leq C\kappa \|x-q_0\|^2$, and hence $|u(x)| \leq C\kappa R^2$ on
	this ball. Since $p=\pi_P p+\ip{p-q_0,\nu}\nu$ for every
	$p\in\RR^{n+1}$, on the cylinder $\pi_P^{-1}(B_{cR}^{P}(q_0))$ the function
	$f(p) = \ip{p-q_0,\nu} - u(\pi_P p)$ is a defining function of the graph of $u$.
	The intrinsic ball $B_{cR}^{\Sigma_1}(p_0)$ lies in this cylinder, because
	$\pi_P p_0 = q_0$ and
	$\|\pi_P p - q_0\| \leq \|p - p_0\| \leq d^{\Sigma_1}(p,p_0)$ for
	$p \in \Sigma_1$.
	Since $\Sigma_1$ is disjoint from $\Sigma_2$ we have that on the connected
	set $B_{cR}^{\Sigma_1} (p_0)$, $f$ is never zero, and hence $f$ is strictly
	positive on that ball, since $f(p_0) = a > 0$.

	Choose the orientation of $\Sigma_1$ so that
	$\ip{\nu_{\Sigma_1}(p_0),\nu}\geq0$ and set
	$\Delta = \|\nu_{\Sigma_1}(p_0) - \nu\|$, which we seek to estimate. We may assume that $\Delta>0$.
	Writing $\beta=\|\pi_{T_{p_0}\Sigma_1}\nu\|$, we have
	\[
		\beta^2
		=1-\ip{\nu_{\Sigma_1}(p_0),\nu}^2
		\geq \frac{1}{2}\Delta^2.
	\]
	Let $\gamma$ be the geodesic starting at $p_0$ with
	$\gamma'(0)=-\beta^{-1}\pi_{T_{p_0}\Sigma_1}\nu$. Properness and the boundary
	separation at $p_0$ allow us to follow $\gamma$ for time $t=cR/2$. Then
	$d^{\Sigma_1}(\gamma(t),p_0)\leq t<cR$, so $\gamma(t)$ lies in
	$B_{cR}^{\Sigma_1}(p_0)$, where $f>0$.
	By \cref{prop-curvature-scale-basic} part $(ii)$ and the bound on $u$,
	\[
		0<f(\gamma(t)) \leq a-\beta t+\frac{\kappa t^2}{2}+C\kappa R^2
		\leq Cr-c\Delta R+C\kappa R^2.
	\]
	Isolating $\Delta$ gives
	\[
		\Delta \leq C r R^{-1} + C \kappa R.
	\]
	Combining this estimate with the normal oscillation estimates centered at
	$p_0$ and $q_0$, for arbitrary $p \in \Sigma_1$ and $q \in \Sigma_2$ we get
	\[
		\|\nu_{\Sigma_1}(p) - \nu_{\Sigma_2}(q)\|
		< C r R^{-1} + C \kappa R. \qedhere
	\]
\end{proof}

%% file: sections/chord_arc_bound.tex
\section{Chord-Arc Bound}\label{sec:chord-arc}
The goal of this section is to prove the chord-arc bound,
\cref{lem:intro-chord-arc}. It says that, at sufficiently low curvature and on the
correct scale, intrinsic distance on a graphical hypersurface is
controlled by ambient distance.

\subsection*{The Euclidean model}
We begin by proving a model Euclidean case, where a stronger result holds. We
will then reduce the general case to the model Euclidean case.

We now describe the model setting. Fix $x = (0,\vec{0}), y = (r,\vec{0})$ where
$r$ is fixed. We denote by $L$ the line segment between $x$ and $y$, and by
$T_\varepsilon$ the tube of radius $\varepsilon$ around $L$. We will assume
that $\varepsilon \leq r$. We can now write our result precisely.

\begin{lemma}[Euclidean Chord-Arc Bound]\label{lem-euclidean-chord-arc}
	There exist constants $\delta_0 > 0$ and $C > 0$ such that the following
	holds. For any graphical hypersurface $\Sigma$ over $T_\varepsilon$, let
	$\kappa > \|A_\Sigma\|_{C^0}$ and define
	the scale-invariant quantity $\delta$ as
	\[
		\delta := r^2 \kappa \varepsilon^{-1}.
	\]
	If $\delta \leq \delta_0$, then
	\[
		\frac
		{d_\Sigma(\hat{x},\hat{y})}
		{d_{T_\varepsilon \times \RR}(\hat{x},\hat{y})}
		\leq 1 + C \delta.
	\]
\end{lemma}
The proof splits according to whether the intrinsic distance is already
short, whether the graph is monotone along the base segment, or whether
one must use level sets to build an efficient path.
\begin{proof}
	By rescaling we may assume that $r = 1$, so $\kappa \varepsilon^{-1} =
		\delta \leq \delta_0$. We may also assume WLOG that $\hat{x}_{n+1} \leq \hat{y}_{n+1}$. We will also write $\rho =
		d_{\Sigma}(\hat{x},\hat{y})$ and $R = d_{T_\varepsilon \times \RR}
		(\hat{x},\hat{y})$. We pick the normal for $\Sigma$ pointing down.

	\paragraph{Case 1.} $\rho < \delta^{-1/2}$.

	In this case we note that we have $\kappa \rho^2 < \varepsilon$ so
	applying \cref{cor-geodesically-connected} with $D = T_\varepsilon
		\times \RR$ shows that $\hat{x}$ and $\hat{y}$ are connected by an
	interior minimizing geodesic on $\Sigma$. Applying
	\cref{prop-curvature-scale-basic} part $(ii)$ gives
	\[
		R \geq \rho\left(1-\frac{\kappa^2\rho^2}{6}\right)
		\geq \rho\left(1-\frac{\delta}{6}\right),
	\]
	since $\kappa^2\rho^2<\kappa^2/\delta=\delta\varepsilon^2\leq\delta$.
	Thus $\rho/R\leq(1-\delta/6)^{-1}\leq 1+C\delta$ for $\delta_0$ small enough.

	\paragraph{Case 2.} $\rho \geq \delta^{-1/2}$. First
	consider the path $\hat{L}$, that is the lift to the surface of the
	straight line between $x$ and $y$. Let $u$ be the function on
	$T_\varepsilon$ whose graph is $\Sigma$, then
	\[
		\rho \leq \|\hat{L}\| = \int_0^1 \sqrt{1 + |\partial_1
			u(t,\vec{0})|^2} d t \leq \max_{L} \sqrt{1 + |\partial_1 u|^2},
	\]
	but we then know that at some point $p$ we have $(1 + \|\nabla u\|^2)^{-1}
		\leq \delta$. Since $T_\varepsilon \seq B(p,2)$ and $\kappa =
	\delta\varepsilon \leq \delta$, \cref{cor-gradient-oscilation} gives
	$\max_{T_\varepsilon} (1 + \|\nabla u\|^2)^{-1} \leq \delta + 4\kappa \leq 5 \delta$. Rearranging
	and again assuming $\delta_0$ is small enough we have that $\|\nabla u\|
		\geq \delta^{-1/2}/3$ everywhere on $T_\varepsilon$.

	\paragraph{Case 2a.} Assume further that $\|\pi_{L \times
			\RR} \nu_\Sigma\| \geq c\varepsilon$ everywhere on $\hat{L}$.
	The path $\hat{L}$ can also be seen as the intersection of the vertical strip
	$L \times \RR$ with $\Sigma$. We have by \cref{prop-curvature-scale-basic} part $(iii)$
	that $\hat{L}$ has curvature as a planar curve bounded above by $C \kappa
		\varepsilon^{-1} < C \delta$. Then applying \cref{cor-gradient-oscilation}
	to the planar curve $\hat{L}$, the graph of $t \mapsto u(t,\vec{0})$,
	we get that since $\max_{L} \sqrt{1 + |\partial_1 u|^2} \geq \delta^{-1/2}$
	then $\sqrt{1 + |\partial_1 u|^2} \geq c\delta^{-1/2}$
	along all of $L$. For $\delta_0$ small enough this implies that
	$\partial_1 u$ is non-zero and hence $u$ is monotonic on $L$. Assume then
	WLOG that it is increasing, then we may apply a $u$ substitution to the
	integral for $\rho$ to get
	\begin{align*}
		\rho \leq & \int_L \sqrt{1 + |\partial_1 u|^2} d z
		= \int_{u(x)}^{u(y)} \sqrt{1 + |\partial_1 u|^2} \frac{d u}{|\partial_1 u|}
		= \int_{u(x)}^{u(y)} \sqrt{1 + |\partial_1 u|^{-2}} d u
		\\ \leq & \int_{u(x)}^{u(y)} \sqrt{1 + C\delta} d u
		\leq |u(y) - u(x)| (1 + C \delta)
		\\ \leq & (1 + C\delta) \sqrt{|u(y) - u(x)|^2 + 1}
		= R (1 + C\delta).
	\end{align*}

	\paragraph{Case 2b.} Assume $\|\pi_{L \times
			\RR} \nu_{\Sigma}\| < c\varepsilon$ at some point $\hat{p} \in
		\hat{L}$. This is the trickiest case, and to tackle it we shall
	need some notation.

	Since we have a lower bound on the gradient, every level set of $u$ is a
	sub-manifold. For a point $z \in T_\varepsilon$, we use $M_z$ to denote the
	connected component of the level set of $u$ at $u(z)$ which passes through
	$z$, with normal chosen in the $u$-increasing direction. The level set is given by $\Sigma \cap \{ x_{n+1} = u(z) \} - u(z)
		e_{n+1}$, it thus has
	\[
		\|\pi_{\{ x_{n+1} = u(z) \}} \nu_\Sigma\| =
		\frac{\|\nabla u\|}{\sqrt{1 + \|\nabla u\|^2}} \geq
		\frac{\delta^{-1/2}/3}{\sqrt{1 + \delta^{-1}/9}} \geq
		\frac{1}{\sqrt{1 + C\delta}}.
	\]
	Applying \cref{prop-curvature-scale-basic} part $(iii)$ we get that $\|A_{M_z}\| <
		\kappa \sqrt{1 + C\delta} \leq \kappa (1 + C\delta/2)$, and by assuming
	$C\delta_0 < 2$ we get that $\|A_{M_z}\| < 2 \kappa$. That is, every level
	set of $u$ has curvature controlled at the relevant scale. Furthermore, at $p$ given in the assumption of
	the case,
	\[
		\|\pi_{L} \nu_{M_p}\| = \frac{\|\pi_{L} \pi_{T_\varepsilon}
			\nu_\Sigma\|}{\|\pi_{T_\varepsilon} \nu_\Sigma\|} \leq
		c\varepsilon.
	\]
	We now carefully consider $M_p$. Restricting to $B_2^{M_p} (p)$ and choosing
	$\delta_0$ small gives $\|\pi_{L} \nu_{M_p}\| \leq c\varepsilon$. Now
	define $r_o$ to be the distance function
	to the $x_1$ axis and note that its gradient is orthogonal to $e_1$, and is
	of magnitude $1$ wherever $r_o \neq 0$. We also define
	\[
		V := \nabla^{M_p} x_1 =
		e_1 - \ip{e_1,\nu_{M_p}}\nu_{M_p}.
	\]
	Then, throughout $B_2^{M_p}(p)$,
	\[
		|V|^2 = V(x_1) = 1 - \ip{e_1,\nu_{M_p}}^2
		\geq 1 - c\varepsilon^2 \geq \frac{3}{4}.
	\]
	On the other hand, wherever $r_o \neq 0$, we have
	\[
		V(r_o)
		= \ip{\nabla r_o,V}
		= -\ip{e_1,\nu_{M_p}}\ip{\nabla r_o,\nu_{M_p}},
	\]
	and hence $|V(r_o)| \leq c \varepsilon$. We can apply this inequality in
	the weak sense across $r_o = 0$. Let $\gamma$ be the integral curve of $V$
	starting at $p$. Since $r_o(p) = 0$, after adjusting $c$ by a fixed
	factor we have $r_o(\gamma(t)) \leq c \varepsilon$ for all $t \in [-2,2]$
	for which $\gamma$ is defined. Thus $\gamma$ does not intersect the sides of
	the $\varepsilon$-tube, so it can only cease to be defined by exiting
	through the caps. Since $d(x_1 \circ \gamma)/dt = V(x_1) \geq 3/4$, for
	every $s \in [-\varepsilon/2, 1+\varepsilon/2]$ there is some $t \in [-2,2]$
	with $x_1(\gamma(t)) = s$ and $r_o(\gamma(t)) \leq c \varepsilon$.

	Now define the slightly shrunk domain $D = [- \varepsilon/2, 1 +
		\varepsilon/2] \times B_{c \varepsilon}^{\RR^{n-1}} (0)$. For any $q \in
		D$, let $s = x_1(q)$. By the construction above, we may choose a point
	$p_q \in M_p \cap B_2^{M_p}(p)$ such that $x_1(p_q) = s$ and, after
	adjusting $c$ by a fixed factor, $\|p_q - q\| \leq c \varepsilon$.
	We compare $M_q$ to $M_p$ locally at $p_q$. If
	$B_{2 \varepsilon}^{M_q} (q)$ and $B_{2 \varepsilon}^{M_p} (p_q)$ intersect,
	then they must be points in the same level set and distance at most
	$4 \varepsilon$ in that level set, so
	$\|\nu_{M_p} (p_q) - \nu_{M_q} (q)\| < C \kappa \varepsilon$. On the other
	hand, if they are disjoint, both local pieces meet
	$B_{c \varepsilon}(q)$. After passing to
	the connected components in $B_{\varepsilon/4}(q)$ which are clearly
	properly embedded, we can apply
	\cref{prop-disjoint-curvature-scale} with $R = \varepsilon/(C 4)$ and
	$r = c \varepsilon$ to the two balls, choosing $c$ small enough so that
	$r$ is below the threshold $cR$ of that proposition, and get that
	$\|\nu_{M_p} (p_q) - \sigma_q\nu_{M_q} (q)\| < C c + C \kappa \varepsilon$
	for a sign $\sigma_q \in \{\pm 1\}$.
	By the normal Lipschitz bound on $M_p$,
	$\|\nu_{M_p}(p_q) - \nu_{M_p}(p)\| \leq C \kappa$, and so, after shrinking
	$\delta_0$ if necessary, in both cases there is a sign $\sigma_q$ with
	\[
		\|\nu_{M_p} (p) - \sigma_q\nu_{M_q} (q)\|
		< C c + C\delta_0.
	\]
	For $c$ and $\delta_0$ small the right hand side is less than $1$, so the
	two sets of $q \in D$ with $\|\nu_{M_q}(q) \mp \nu_{M_p}(p)\| < C c +
	C\delta_0$ are disjoint. They are open, since $q \mapsto \nu_{M_q}(q) =
	\nabla u(q)/\|\nabla u(q)\|$ is continuous, they cover the connected set
	$D$, and the first one contains $p$. Hence $\sigma_q = 1$ for all $q \in D$, that is,
	$\|\nu_{M_p}(p) - \nu_{M_q}(q)\| < Cc + C\delta_0$ on $D$.

	Define the fixed unit vector $w := \frac{\pi_{L^\perp} \nu_{M_p} (p)}{\|\pi_{L^\perp} \nu_{M_p} (p)\|}$
	then using \cref{prop-simple-geometric} part $(iv)$ we get
	\[
		\|\nu_{M_p} (p) - w \| \leq \sqrt{2} |\ip{e_1, \nu_{M_p}}| \leq
		C c \varepsilon.
	\]
	Another application of \cref{prop-simple-geometric} part $(iv)$ gives
	\[
		\|\nu_{M_q}(q) - \nu_{\Sigma}(q)\|
		\leq \sqrt{2} |\ip{\nu_{\Sigma}(q), e_{n+1}}|
		\leq \sqrt{10 \delta}.
	\]
	Combining all of the above, and using $\delta = \kappa \varepsilon^{-1}$ and
	$\varepsilon \leq 1$, we have
	\begin{align*}
		\|\nu_{\Sigma}(q) - w\|
		 & \leq C\sqrt{\delta} + C c + C\delta
		\\
		 & \leq C\sqrt{\delta_0} + C c + C\delta_0.
	\end{align*}
	Choose $c$ and $\delta_0$ sufficiently small so that the last quantity is
	small. We conclude that
	$\|\nu_{\Sigma} (q) - w\| \leq c$ for all $q \in D$. Now
	$T_{c \varepsilon} \seq D$ and so we
	may freely assume that $\|\nu_{\Sigma} (q) - w\| \leq c$ everywhere in
	$T_\varepsilon$ after replacing $\varepsilon$ with $c\varepsilon$ and
	$\delta_0$ with $c\delta_0$.

	Here the constants are chosen in this order: first the radius
	$c_b\varepsilon$ of $D$, then the threshold $c_c\varepsilon$ for
	$\|\pi_L\nu\|$ in the claim below, small enough that the flow argument
	above works from a point of $L$ where $\|\pi_L\nu\| < 3c_c\varepsilon$,
	then the constant $c_a := c_bc_c/2$ separating Cases 2a and 2b, so that $p$
	still satisfies $\|\pi_L\nu_{M_p}(p)\| \leq 2c_a\varepsilon = c_c
	(c_b\varepsilon)$ in the shrunk tube, and finally $\delta_0$.

	In particular, $\|\nabla_w u\| \geq c \|\nabla u\|$, and so
	$u$ is strictly increasing in the direction of $w$.

	We now consider $\Sigma' = \Sigma \cap \Span(e_1, w, e_{n+1})$, by the
	above we get that $\|\pi_{\Span(e_1,w,e_{n+1})} \nu\| \geq 1 - c \geq
		\frac{1}{2}$. We thus have by \cref{prop-curvature-scale-basic} that $\|A_{\Sigma'}\|
		\leq 2 \kappa$ and so, after replacing $\delta_0$
	with $\delta_0/2$ we may assume that $n = 2$. In the $n = 2$ case we have
	$e_2 = w$ and $e_3 = e_{n+1}$ and we will use this coordinate system moving
	forward.

	A direct consequence of this simpler setup is that since $n = 2$ then
	the level sets of $u$ are curves. Since $\|\nu_\Sigma - e_2\| \leq c$,
	their arclength parametrizations $\gamma(t)$, oriented in the increasing
	$x_1$ direction, satisfy $\partial_t x_1(\gamma(t)) \geq 1/2$ for $c$
	small enough. Thus every level-set segment between $x_1 = 0$ and
	$x_1 = 1$ has length at most $2$.

	\begin{claim}
		There exists a point $q \in L$ such that $u(x) \leq u(q) \leq
			u(y)$, the level set $M_q$ intersects every line $x_1 = s$ for $s
			\in [0, 1]$, and $\|\pi_{L} \nu_{M_q}\| < c \varepsilon$ at
		some point $z$ in $M_q$ with $z_1 \in [0,1]$.
	\end{claim}
	\begin{proof}
		Assume this is false. We already know that $p$, which we had
		assumed to exist earlier, satisfies the second condition. It must
		then fail the first, that is either $u(p) < u(x)$ or $u(p) >
			u(y)$. By replacing $u$ with $-u$ and $x$ with $y$ if necessary
		we may assume that $u(p) > u(y)$.

		We thus know that $u$ is not monotonic along $L$, so set $p' \in L$ to
		be the maximum of $u$ along $L$ and note that $\pi_L \nu_{M_{p'}} = 0$.
		Applying the same $V = \nabla^{M_{p'}} x_1$ flow argument as above, we
		know that $M_{p'}$ intersects every $x_1 = s$ for
		$s \in [-\varepsilon/2, 1 + \varepsilon/2]$ while staying within the
		corresponding radial bound, and in particular $|r_o (M_{p'})| \leq
			C \kappa$. Since $u$ is strictly increasing in the
		direction of $e_2$, and since $u(p') \geq u(q)$ for all $q \in L$ we
		know that $M_{p'}$ lies entirely on one side of $L$, in the direction of
		$e_2$.

		Let us restrict $M_{p'}$ to $[0,1] \times [-\varepsilon,\varepsilon]$,
		and then consider the flow $\Phi_t$ of $M_{p'}$ along the vector field
		$-e_2 \|\nabla_{2} u\|^{-1}$. One can easily check that
		$u(\Phi_t(M_{p'})) = u(p') - t$ whenever the flow is well defined. We
		now prove that it exists at least until $t = u(p') - u(y)$. To show
		this, it is enough to show that for all $t \in [0, u(p') - u(y)), q \in
			M_{p'}$, $x_2(\Phi_t (q)) > -\varepsilon/2$.

		We will now construct a barrier that the flow cannot cross for the
		required times, which will prove what is required. This barrier will be
		the short circular arc $\mathcal{C}$ between $x$ and $y$ with radius $(2
			\kappa)^{-1}$, see \cref{fig:geometric-setup}.

		\begin{figure}[htbp]
			\centering
			\input{figures/geometric-setup.tex}
			\caption{Geometric setup for when $u$ is not monotonic along $L$}
			\label{fig:geometric-setup}
		\end{figure}
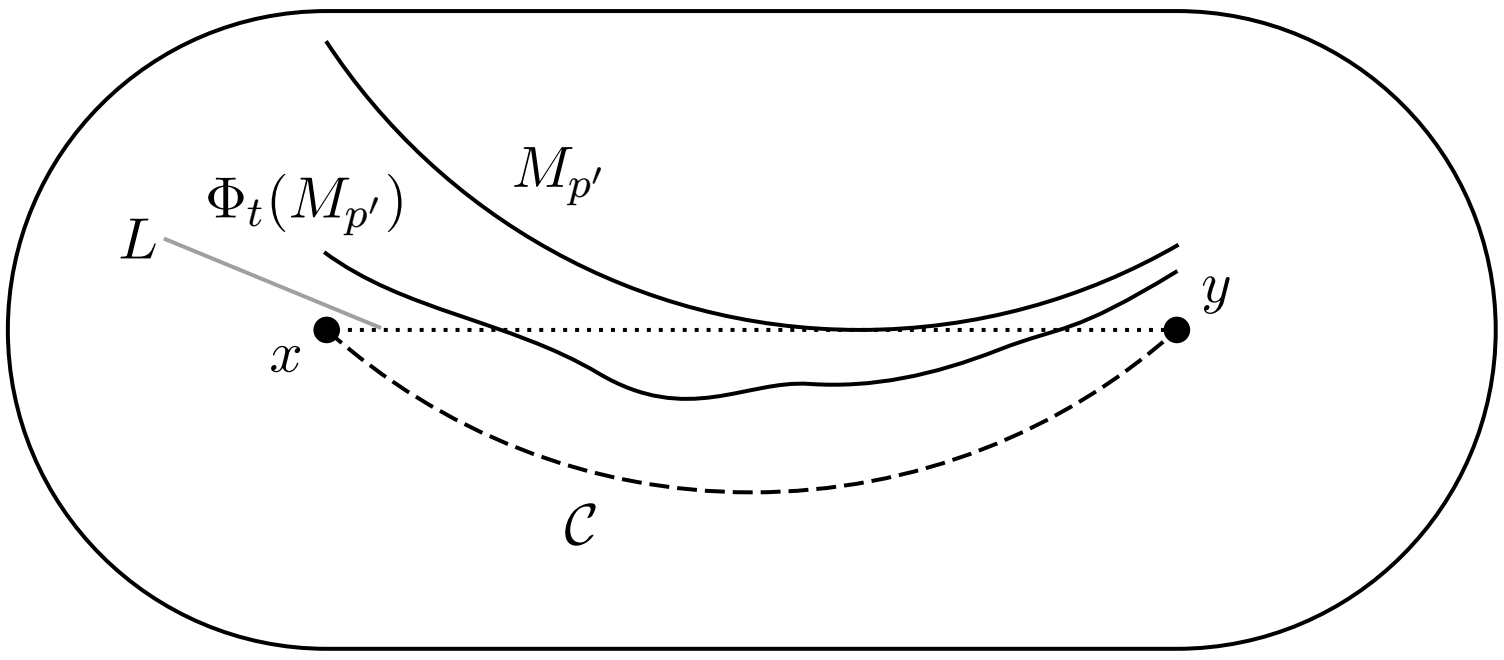
		To see that $\Phi_t (M_{p'})$ cannot intersect $\mathcal{C}$, assume for a
		contradiction that $t_*$ is the first time for which it does.
		Since $t_* < u(p') - u(y)$ we know that $u(\Phi_{t_*} (M_{p'})) =
			u(p') - t_* > u(y)$ and hence the intersection cannot happen at
		$x$ or at $y$. Thus, the intersection happens at the interior at
		some point $q$. Since $t_*$ is the first $t$ where this happens,
		we must have that $\mathcal{C}$ is tangent to $\Phi_{t_*} (M_{p'})$ at $q$.

		Now $\mathcal{C}$ has curvature $2 \kappa$ and since $\Phi_{t_*} (M_{p'})$ is a
		level set, it has curvature strictly less than $2 \kappa$, so by the
		comparison principle \cref{prop-onesided}, $\Phi_{t_*} (M_{p'})$ cannot
		lie on the top side of $\mathcal{C}$. This is a contradiction since we assumed
		$t_*$ is the first time at which the two intersect. We thus know that
		$\Phi_{t} (M_{p'})$ never intersects $\mathcal{C}$ for $t \in [0,u(p') - u(y))$
		and thus $\Phi_{u(p') - u(y)} (M_{p'}) = M_{y}$ is a level set which is
		bounded between $\mathcal{C}$ and $M_{p'}$ and thus indeed intersects every line
		$x_1 = s$ for $s \in [0,1]$. Since $u(x) \leq u(y) \leq u(p')$ we know
		that $M_y$ dips below $L$ where $L$ intersects $M_{p'}$ and is not
		below $L$ at $y$ and at $x_1 = 0$. Thus we know that $x_2$ achieves an
		interior minimum on $M_y$ and at that point we have $\nu_{M_y} = e_2$ so
		$\pi_L \nu_{M_y} = 0$. So $y$ is a $q$ for our claim.
	\end{proof}
	With the claim established we assume that $q$ is a point as
	established in the claim. First we can construct a relatively
	efficient path between $\hat{x}$ and $\hat{y}$. We will define the
	path $\gamma$ in $T_\varepsilon$ and then lift to a path
	$\hat{\gamma}$ in $\Sigma$.

	First denote by $q_1,q_2$, the points on $M_q$ with $x_1 = 0$ and $x_1 = 1$
	respectively. If there are several pick the ones closest to $x$ and $y$
	respectively. Now define the following 3-piece path in $T_\varepsilon$,
	take the straight line $x \to q_1$, which we will denote $\gamma_1$. Next
	take the shortest path in $M_q$ between $q_1$ and $q_2$ which we will
	denote $\gamma_2$. Finally the straight line $q_2 \to y$ we denote
	$\gamma_3$. We then define $\gamma$ to be their concatenation.

	We now can upper bound the length of $\hat{\gamma}$ by upper bounding the
	lengths of each $\hat{\gamma}_i$. First we can upper bound the lengths of
	$\hat{\gamma}_1,\hat{\gamma}_3$. On each slice of $x_1 = s$, $u$ is
	monotonic, hence we have
	\begin{align*}
		L(\hat{\gamma}_1)
		={} & \int_{0}^{x_2 (q_1)} \sqrt{1 + |\nabla_2 u|_{(0,t)}|^2} d t
		\leq \int_{0}^{x_2 (q_1)} 1 + |\nabla_2 u|_{(0,t)}| d t
		\\={} & \int_{0}^{x_2 (q_1)} 1 + \nabla_2 u|_{(0,t)} d t
		= x_2 (q_1) + u(q) - u(x),
	\end{align*}
	in an identical fashion we also get
	\[
		L(\hat{\gamma}_3) \leq -x_2 (q_2) + u(y) - u(q),
	\]
	and finally $\hat{\gamma}_2 = \gamma_2 + u(q)e_{3}$ and so the two curves
	have the same length, which is at most $2$ by the level-set segment
	bound above.

	Combining these we get
	\begin{align*}
		L(\hat{\gamma}) \leq{} & x_2 (q_1) + u(q) - u(x) + 2 -x_2 (q_2) +
		u(y) - u(q) \leq u(y) - u(x) + 2 + 2\varepsilon
		\\ \leq {} & u(y) - u(x) + 4.
	\end{align*}
	This gives us a lower bound on $u(y) - u(x)$, via $u(y) - u(x) \geq
		L(\hat{\gamma}) - 4 \geq \rho - 4$. For brevity we will denote $h =
		|u(y) - u(x)|$.

	We can now estimate the distance quotient using this path,
	\[
		\left(\frac{\rho}{d_{T_\varepsilon \times \RR} (\hat{x},\hat{y})}\right)^2
		= \frac{\rho^2}{1 + h^2}
		\leq \frac{h^2 + 8 h + 16}{1+h^2}
		= 1 + \frac{8h + 15}{1+h^2} \leq 1 + \frac{C}{h} \leq 1 + \frac{C}{\rho}.
	\]
	Then since $\sqrt{1 + x} \leq 1 + x/2$ we get that the quotient is bounded
	by $1 + C \rho^{-1}$, and so since we assumed $\rho \geq \delta^{-1/2}$ we
	know that the quotient is bounded by $1 + C \sqrt{\delta}$. Now if $\rho
		\geq c \delta^{-1}$, then in fact the ratio is bounded by $1 + C
		\delta$ and so this path is efficient enough. We hence restrict
	ourselves to the case $\delta^{-1/2} \leq \rho \leq c \delta^{-1}$.

	By the claim, some point $z$ in $M_q$ has $\|\pi_L \nu_{M_q}\| \leq c
		\varepsilon$, and hence the same is true for $\|\pi_L \nu_{\Sigma}\|$ at
	$\hat{z}$. Now via the same construction as for the path $\gamma$, for any
	point $\hat{a} \in \Sigma$ with $x_1(\hat{a}) \in [0,1]$ and $x_3(\hat{a})
		\in [u(x),u(y)]$ we can construct a path between $\hat{z}$ and $\hat{a}$
	with length at most $h + 2 + 2\varepsilon$, which is bounded above by $2 \rho$. Thus we
	get $\|\nu_\Sigma (\hat{a}) - \nu_\Sigma (\hat{z})\| < 2 \kappa \rho \leq
		c \varepsilon$. Thus we have that
	$\|\pi_{L} \nu_{\Sigma}(a)\| \leq c \varepsilon$. In particular
	this is true for $\hat{x}$ and $\hat{y}$.

	Now using the same reasoning as we did earlier for the level set $M_p$, the
	level sets $M_x$ and $M_y$ also intersect every $x_1 = s$ for $s \in
		[-\varepsilon/2, 1 + \varepsilon/2]$. We now project $\Sigma$ onto the
	strip $L \times \RR$. Since we know that $u$ is strictly monotonic in the
	$x_2$ direction, this map is injective. Since both level sets $M_x,M_y$
	intersect every slice $x_1 = s$ for $s \in [0,1]$, the image of this map
	contains the rectangle $L \times [u(x), u(y)]$. We thus can consider the
	preimage $\pi_{L \times [u(x),u(y)]}^{-1} \Sigma$ as a graph of a
	\textit{different} function $\mu$ over this rectangle. We now consider the
	straight line unit time path between $\hat{x}$ and $\hat{y}$, which we will
	call $\sigma$. This line is actually inside the rectangle $L \times
		[u(x),u(y)]$, thus we can lift it to $\Sigma$ via the function $\mu$, we
	will denote by $\mu(\sigma)$ the lift of $\sigma$ to $\Sigma$.

	To compute the length of $\mu(\sigma)$ we will need to get a good handle on
	$\nabla \mu$ and that can be done via a good handle on $\nu_\Sigma$. We can
	write $\nu_\Sigma = b_1 e_1 + b_2 e_2 + b_3 e_3$. By earlier estimates we
	know that $|b_1| \leq c \varepsilon$, from the start of case $2b$
	we also know that $|b_3| \leq C\sqrt{\delta}$. For $\delta_0$ small enough
	this then implies that $|b_2| \geq 1/2$.

	We then have
	\[
		L(\mu(\sigma)) \leq \int_\sigma \sqrt{\|\partial_t \sigma\|^2 +
			|\nabla_{\partial_t \sigma} \mu|^2} dt.
	\]
	First note that $\|\partial_t \sigma\| = R$, and using our estimates above
	we have
	\[
		\|\nabla_{\partial_t \sigma} \mu\|
		= \|\nabla_1 \mu + h \nabla_3 \mu\|
		\leq \frac{|\ip{e_1, \nu_\Sigma}|}{|\ip{e_2, \nu_\Sigma}|}
		+ h \frac{|\ip{e_3, \nu_\Sigma}|}{|\ip{e_2, \nu_\Sigma}|}
		\leq c \varepsilon + C h \sqrt{\delta}
		\leq C R \sqrt{\delta},
	\]
	hence we get
	\[
		L(\mu(\sigma)) \leq \int_\sigma \sqrt{R^2 + C R^2 \delta} dt
		= R \sqrt{1 + C \delta} \leq R (1 + C \delta),
	\]
	and this finishes the proof.
\end{proof}
\subsection*{Reduction to the Euclidean setting}
With the Euclidean case proven, we reduce the Riemannian estimate to this
model setting in a small metric tube. Here we write $r$ to mean $d_N (x,y)$.
\begin{theorem}[Chord-Arc Bound, \cref{lem:intro-chord-arc}]\label{lem-chord-arc}
	There exist constants $\delta_0 > 0$ and $C > 0$ such that the following
	holds. With $\kappa,K$ as in \eqref{eq-kappa}
	and \eqref{eq-ambient-curvature}, for all $x,y \in \Omega$ and every
	minimizing geodesic $\gamma$ of $N$ from $x$ to $y$, if
	$0 < \varepsilon \leq r$ and $B_\varepsilon(\gamma) \seq \Omega$, define the scale-invariant
	quantity $\delta$ as
	\[
		\delta := r^2 \left(\kappa\varepsilon^{-1} + K\right).
	\]
	If $\delta \leq \delta_0$, then $\Sigma$ satisfies the Chord-Arc estimate
	\[
		\frac{d_{\Sigma} (\hat{x}, \hat{y})}{d_{N \times \RR} (\hat{x},
			\hat{y})} \leq 1 + C \delta.
	\]
\end{theorem}
\begin{proof}
	We prove by reduction to the Euclidean case.

	Extend $\gamma$ by a length of $\varepsilon/2$ at both end points. Every
	point within distance $\varepsilon/2$ of the extended geodesic lies
	within distance $\varepsilon$ of $\gamma$, hence in $\Omega$.

	Since $\varepsilon^2 K \leq r^2K \leq \delta_0$, applying \cref{prop-connection-expansion} at scale $\varepsilon/2$ to the extended geodesic gives a map $\Phi : U_{c \varepsilon/2} \to \Omega$ and a pullback metric $\tilde{g}$ which is $C^1$ close to the flat metric $\ov{g}$. We now pullback the function $u$ to $U_{c \varepsilon/2}$, to give a graphical hypersurface $\ov{\Sigma}$ over $U_{c \varepsilon/2}$.
	Set $\ell=d_N(x,y)$ and denote the corresponding points of $\ov{\Sigma}$
	by $\tilde{x}=(0,u(x))$ and $\tilde{y}=(\ell e_1,u(y))$.
	The estimates in \cref{prop-connection-expansion} then give
	\[
		\|\tilde{g} - \ov{g}\|_{C^0} \leq C \varepsilon^2 K,\qquad
		\|\nabla^{\tilde{g}} - \ov{\nabla}\|_{C^0} \leq C \varepsilon K.
	\]
	Combining these with \cref{prop-h-bounds} gives
	\[
		\|\ov{h}_{\ov{\Sigma}}\|_{C^0}
		\leq C(\kappa + \varepsilon K)
		\leq C \delta \varepsilon r^{-2}.
	\]
	With this bound now in flat coordinates, we find that
	$r^2\|\ov{h}_{\ov{\Sigma}}\|_{C^0}\bar{\varepsilon}^{-1}
		\leq C \delta$ where $\bar{\varepsilon} = c \varepsilon/2$. After decreasing $\delta_0$ so this quantity is below
	the threshold in \cref{lem-euclidean-chord-arc}, the Euclidean estimate gives
	\[
		\frac{\ov{d}_{\ov{\Sigma}}(\tilde{x},\tilde{y})}{\ov{d}_{T_{\bar{\varepsilon}} \times \RR}
			(\tilde{x}, \tilde{y})} \leq 1 + C \delta.
	\]
	Let $\sigma$ be a path in $\ov{\Sigma}$ with length
	$\ov{L}(\sigma)\leq(1+C\delta)\ov{d}_{T_{\bar{\varepsilon}}\times\RR}
		(\tilde{x},\tilde{y})$. Then $\Phi \times\Id_{\RR} (\sigma)$ gives a path in
	$\Sigma$ from $\hat{x}$ to $\hat{y}$. Then we have
	\begin{align*}
		L_g(\Phi \times\Id_{\RR} (\sigma)) = {} & L_{\tilde{g}} (\sigma) \leq (1 + \|\tilde{g} - \ov{g}\|_{C^0}) L_{\ov{g}}(\sigma) \leq (1 + C \varepsilon^2 K) (1 + C \delta) \ov{d}_{T_{\bar{\varepsilon}}\times\RR} (\tilde{x},\tilde{y})
		\\
		= {}                                    & (1 + C \varepsilon^2 K) (1 + C \delta) d_{N \times\RR} (\hat{x},\hat{y}) \leq (1 + C \delta) d_{N \times\RR} (\hat{x},\hat{y}),
	\end{align*}
	where in the last inequality we used that $\varepsilon^2 K \leq r^2 K \leq \delta$.

	We thus get
	\[
		d_\Sigma (\hat{x}, \hat{y})
		\leq (1 + C \delta)d_{N \times \RR} (\hat{x},\hat{y}),
	\]
	hence the theorem is proved for some small enough choice of
	$\delta_0$ and large enough choice of $C$.
\end{proof}

%% file: figures/geometric-setup.tex
\tikzset{every picture/.style={line width=0.75pt}} 
\begin{tikzpicture}[x=1pt,y=1pt,yscale=1,xscale=1,cm={0,1,-1,0,(0,0)}]

	\draw   (70,130) .. controls (70,96.86) and (96.86,70) ..
	(130,70) .. controls (163.14,70) and (190,96.86) .. (190,130)
	-- (190,290) .. controls (190,323.14) and (163.14,350) ..
	(130,350) .. controls (96.86,350) and (70,323.14) .. (70,290) -- cycle ;
	\draw  [draw opacity=0][fill={rgb, 255:red, 0; green, 0;
				blue, 0 }  ,fill opacity=1 ] (127.5,290) .. controls
	(127.5,288.62) and (128.62,287.5) .. (130,287.5) .. controls
	(131.38,287.5) and (132.5,288.62) .. (132.5,290) .. controls
	(132.5,291.38) and (131.38,292.5) .. (130,292.5) .. controls
	(128.62,292.5) and (127.5,291.38) .. (127.5,290) -- cycle ;
	\draw  [dash pattern={on 3.75pt off 1.5pt}] (130,290) ..
	controls (111,268.77) and (99.44,240.73) .. (99.44,210) ..
	controls (99.44,179.27) and (111,151.23) .. (130,130) ;
	\draw [color={rgb, 255:red, 0; green, 0; blue, 0 }  ,draw
		opacity=1 ] [dash pattern={on 0.75pt off 1.5pt}]  (130,130)
	-- (130,290) ;
	\draw  [draw opacity=0][fill={rgb, 255:red, 0; green, 0;
				blue, 0 }  ,fill opacity=1 ] (127.5,130) .. controls
	(127.5,128.62) and (128.62,127.5) .. (130,127.5) .. controls
	(131.38,127.5) and (132.5,128.62) .. (132.5,130) .. controls
	(132.5,131.38) and (131.38,132.5) .. (130,132.5) .. controls
	(128.62,132.5) and (127.5,131.38) .. (127.5,130) -- cycle ;
	\draw    (141.09,129.9) .. controls (128.69,150.21) and
	(131.04,151.25) .. (126.24,163.57) .. controls
	(121.44,175.89) and (119.04,187.25) .. (119.84,199.17) ..
	controls (120.64,211.09) and (111.84,222.02) ..
	(121.84,238.85) .. controls (131.84,255.67) and
	(133.44,275.43) .. (144.64,290.47) ;
	\draw [color={rgb, 255:red, 160; green, 160; blue, 160 }
		,draw opacity=1 ]   (130.35,279.75) -- (147.2,320.65) ;
	\draw   (184.34,290.12) .. controls (151.62,268.69) and
	(130,231.7) .. (130,189.67) .. controls (130,167.82) and
	(135.84,147.34) .. (146.04,129.69) ;

	\draw (128,293.4) node [anchor=north east] [inner sep=0.75pt]    {$x$};
	\draw (132,126.57) node [anchor=south west] [inner sep=0.75pt]    {$y$};
	\draw (147.2,320.65) node [anchor=east] [inner sep=0.75pt] {$L$};
	\draw (98.27,237.71) node [anchor=north east] [inner
		sep=0.75pt]    {$\mathcal{C}$};
	\draw (152.4,256.3) node [anchor=south west] [inner
		sep=0.75pt]    {$M_{p'}$};
	\draw (146.64,293.87) node [anchor=south][inner sep=0.75pt]
	{$\Phi _{t}( M_{p'})$};
\end{tikzpicture}

%% file: sections/inscribed_curvature.tex
\section{In-Radius Bound}\label{sec:main}
Let us rigorously define the in-radius. We first study the properties of
tangent balls to make sure the in-radius is well defined for an arbitrary
hypersurface.
\begin{definition}\label{def-tangent-ball}
	Let $\Sigma \seq N$ be any oriented hypersurface in a complete
	Riemannian manifold. Let $p \in \Sigma$ be an arbitrary interior
	point and $\nu(p)$ the normal vector at $p$. For any $r > 0$, the
	interior tangent ball of radius $r$ at $p$ is given by the (open)
	ball $B(\exp_p(-r \nu(p)), r)$ which we denote $TB_p (r)$.
\end{definition}
\begin{proposition}\label{prop-tangent-ball-monotonicity}
	If $r_1 < r_2$ then $TB_p (r_1) \seq TB_p (r_2)$.
\end{proposition}
\begin{proof}
	Let $x \in TB_p (r_1)$ then $d(x, \exp_p(-r_1 \nu(p))) < r_1$, but then
	\begin{align}
		d(x, \exp_p(-r_2 \nu(p))) & \leq d(x, \exp_p(-r_1 \nu(p))) +
		d(\exp_p(-r_1 \nu(p)), \exp_p(-r_2 \nu(p)))
		\\                     & < r_1 + (r_2 - r_1) = r_2,
	\end{align}
	hence $x \in TB_p (r_2)$.
\end{proof}

Now we want our balls to not just be tangent but to be entirely on one side of
the hypersurface. We thus define the following notion.
\begin{definition}\label{def-admissible}
	An interior tangent ball $TB_p (r)$ is \textbf{admissible} if $TB_p
		(r) \cap \Sigma = \varnothing$.
\end{definition}
From the monotonicity property above we can deduce that if $TB_p (r)$
is admissible then $TB_p (r')$ is also admissible for all $r' < r$.
Similarly if $TB_p (r)$ is not admissible then $TB_p (r')$ is not
admissible for all $r' > r$.

\begin{proposition}\label{prop-admissible-existence}
	Let $\Sigma \seq N$ be any oriented embedded hypersurface,
	with or without boundary, in a complete Riemannian manifold $N$. Let $p$ be an
	interior point of $\Sigma$, there exists an $r$ such that $TB_p
		(r)$ is admissible.
\end{proposition}
\begin{proof}
	Since $\Sigma$ is embedded, it is properly embedded in some open
	neighborhood $W$ of $p$. By the tubular neighborhood theorem for properly
	embedded hypersurfaces, a
	neighborhood of $p$ in the normal bundle embeds into $W$ and has nearest-point
	projection to $\Sigma \cap W$. Let $r>0$ be small enough so that
	$B(p,3r) \seq W$ and
	$q := \exp_p(-r\nu(p))$ lies in this tubular neighborhood and projects to
	$p$. Then $d(q,\Sigma)=r$. Thus for every $x \in \Sigma$ we have
	$d(q,x)\geq r$, and since $TB_p(r)=B(q,r)$ is open, it follows that
	$TB_p(r)\cap\Sigma=\varnothing$.
\end{proof}

The picture is now clearer. All sufficiently small tangent balls are
admissible, they either remain admissible for all $r$ or at some critical $r$
they stop being admissible. This motivates the following definition.
\begin{definition}\label{def-in-radius}
	The in-radius of $\Sigma$ at $p$ is the extended real number $r_* \in
		(0,\infty]$ given by
	\[
		r_* = \sup \{r : TB_p (r)\, \text{is admissible}\}.
	\]
\end{definition}

The next lemma gives the precise alternative governing the in-radius at a point
$p$: either the cut locus or another contact point controls the admissible
tangent balls, or $\partial TB_p(r_*)$ agrees with $\Sigma$ to second order
at $p$ in some direction.
\begin{lemma}[In-Radius
		Trichotomy]\label{lem-inscribed-curvature-trichotomy}
	Let $\Sigma$ be an oriented hypersurface embedded in a
	complete manifold $N$. Let $p$ be an interior point to $\Sigma$. If $r_* <
	\infty$, assume that $\Sigma \cap \ov{TB_p(r)}$ is compact for some $r >
	r_*$; this holds, for instance, if $\Sigma$ is properly embedded.
	Then at least one of the following holds:
	\begin{enumerate}
		\item $r_* = t_{cut} (p, -\nu(p))$.
		\item $r_* < t_{cut} (p, -\nu(p))$, and there exists a point $q
			      \neq p$ in $\Sigma$ such that $q \in \partial TB_p (r_*)$.
		\item $r_* < t_{cut} (p, -\nu(p))$, and
		      $h_{\partial TB_p(r_*)}(p) (v,v) = h_{\Sigma} (p) (v,v)$ for some
		      unit vector $v \in T_p \Sigma$.
	\end{enumerate}
\end{lemma}
\begin{proof}
	We clearly cannot have $r_* > t_{cut}(p, -\nu(p))$. Assume then
	that $r_* < t_{cut}(p, -\nu(p))$, then we have that $\partial TB_p
		(r)$ is smooth in some neighborhood of $p$ for all
	$r \leq r_* + \varepsilon$, for $\varepsilon$ sufficiently small.
	Now define the set $S =
		\bigcap_{r_* < r < r_* + \varepsilon} \ov{TB_p(r) \cap \Sigma}$.
	By completeness of $N$ and our assumption, after decreasing $\varepsilon$,
	these are closed
	subsets of the compact set $\ov{TB_p(r_*+\varepsilon)} \cap \Sigma$.
	They are nested by \cref{prop-tangent-ball-monotonicity} and non-empty by
	definition of the in-radius, so their total intersection is non-empty.

	Now pick $q \in S$, by definition we can pick for any $r > r_*$, points in
	$\Sigma \cap TB_p(r)$ arbitrary close to $q$. Thus we have $d(q, \exp_p (-r
		\nu(p))) \leq r$ so by continuity we have $d(q, \exp_{p}(- r_* \nu(p))) \leq
		r_*$. On the other hand if $d(q, \exp_p(-r_* \nu(p))) < r_*$ then again by
	continuity for some $r$ slightly smaller than $r_*$ we will have $d(q,
		\exp_p (-r \nu(p))) < r$ and so $q \in TB_p (r)$ which contradicts the fact
	that $r_*$ is the in-radius. We thus have $d(q, \exp_{p} (-r_* \nu(p))) =
		r_*$ and so $q \in \partial TB_p (r_*)$.

	If we can pick $q \neq p$, we are in case $(2)$, so we assume instead that
	$S = \{p\}$. For any $r_* < r < r_* + \varepsilon$, take a sequence
	$x_n \to p$ in $T B_p (r) \cap \Sigma$. These points lie on the side of
	$\partial T B_p (r)$ opposite to $\nu(p)$, so they witness that
	$\partial T B_p (r)$ is not locally interior to
	$\Sigma$ at $p$. Hence the strict converse direction of \cref{prop-onesided}
	implies that
	$h_{\partial T B_p (r)} - h_{\Sigma}$ must have a non-positive eigenvalue.
	Since $\partial TB_p(r)$ is smooth near $p$ for $r \leq r_* + \varepsilon$, its
	second fundamental form at $p$ is continuous in $r$, so
	$h_{\partial TB_p(r_*)}(p) - h_{\Sigma}(p)$ has a non-positive eigenvalue.
	On the other hand, for $0<r<r_*$ the ball $TB_p(r)$ is admissible, so
	$\Sigma$ lies outside it and \cref{prop-onesided} gives
	$h_{\partial TB_p (r)}(p) \geq h_\Sigma(p)$. Letting $r \nearrow r_*$ we
	get $h_{\partial TB_p(r_*)}(p) - h_{\Sigma}(p) \geq 0$, so it has a $0$
	eigenvalue which completes the proof.
\end{proof}

%% file: sections/main_theorem.tex
Using \cref{lem-chord-arc} we can now prove \cref{thm:main}.
\begin{theorem}[In-Radius
		Estimates, \cref{thm:main}]\label{thrm-inscribed-curvature-estim}
	With $\rho,\kappa,K$ as in \eqref{eq-rho}, \eqref{eq-kappa}, and
	\eqref{eq-ambient-curvature}, there exists a dimensional constant $C > 0$
	such that for either orientation of $\Sigma$ and all interior points $p \in
	\Sigma$ we have
	\begin{align*}
		r_*(p)^{-1} \leq C (\kappa + K^{1/2}) + 4 \rho(p)^{-1}.
	\end{align*}
\end{theorem}
\begin{proof}
	Assume for a contradiction that the result is false, that is for
	any $\alpha > 0$ we have for some $(N,\Omega,\Sigma)$ the inequality
	\[
		r_*(p)^{-1} > \alpha^{-1} (\kappa + K^{1/2})
		+ 4 \rho(p)^{-1}.
	\]
	Clearly $r_*(p)$ is finite so via rescaling, we assume
	that $r_*(p) = 1$ at some point $p \in \Sigma$, and hence
	\[
		\alpha > (\kappa + K^{1/2}) \quad
		\text{and} \quad \rho(p) > 4.
	\]
	This implies in particular that $K \leq \alpha^2$. After decreasing
	$\alpha$ we may also assume $\alpha<1$. Since $\rho(p)^{-1} \leq 1/4$ we have
	$\inj(N) \geq 4$, so $r_* = 1 < \inj(N)$.

	We now apply \cref{lem-inscribed-curvature-trichotomy} to $p$. Its
	compactness assumption holds: since $\rho(p) > 4$, the closed ball
	$\ov{B(\pi_N(p),3.5)}$ lies in $\Omega$, so the part of $\Sigma$ over it is
	compact, and for $r \leq 1.75$ the ball $TB_p(r)$ lies over
	$B(\pi_N(p),3.5)$. The first
	alternative is impossible because $r_*<\inj(N)$. Suppose the third
	alternative holds. Then for some unit vector $v \in T_p\Sigma$,
	\[
		h_{\partial TB_p(1)}(p)(v,v)=h_\Sigma(p)(v,v),
	\]
	and therefore $h_{\partial TB_p(1)}(p)(v,v) \leq \kappa \leq \alpha$. On the
	other hand, the sectional curvature of $N \times \RR$ is at most
	$K \leq \alpha^2 < 1$, so the Hessian comparison theorem gives
	\[
		h_{\partial TB_p(1)}(p) \geq \sqrt{K}\cot\bigl(\sqrt{K}\bigr) g
		\geq (1-\alpha^2) g .
	\]
	For $\alpha$ sufficiently small this is a contradiction. Thus the second
	alternative holds: there is a point $q \neq p$ such that a ball of radius
	$r_*$ is tangent to $\Sigma$ at both $p$ and $q$. We denote this ball by
	$B$ and its center by $o$. Note that the normals of $\Sigma$ at $p$ and $q$ agree with the outward normals of $\partial B$ because $\Sigma$ is a graph and so $B$ lies either above or below it.

	The points $p$ and $q$ are at most distance $3$ away from each other so
	both are distance at least $1$ from $\partial \Omega$. Since $\Sigma$ is
	graphical, $p$ and $q$ are the graph lifts over
	$x:=\pi_N(p)$ and $y:=\pi_N(q)$, respectively. Now pick $r = d(x,y)$, and
	then $\varepsilon = r/3$. Then $r>0$, and for a minimizing geodesic
	$\beta$ of $N$ from $x$ to $y$ we have $B_\varepsilon(\beta) \seq
	B(x,4) \seq \Omega$, since $r \leq 3$. So the hypotheses of
	\cref{lem-chord-arc} apply once $\delta$ is below $\delta_0$. We check this, 
	\[
		\delta = r^2 (\kappa \varepsilon^{-1} + K) =
		r^2 (3 \kappa r^{-1} + K) \leq
		9 \alpha + 9 \alpha^2,
	\]
	so by picking $\alpha$ small enough we may assume that $9 \alpha + 9
		\alpha^2$ is less than 
	$\delta_0$ from \cref{lem-chord-arc}. By the result of that theorem we then get
	that $d_{\Sigma}(p,q) \leq (1 + C \alpha) d_{N
				\times \RR}(p,q) $. Assume that $\alpha$ is small enough so that $1 + C
		\alpha \leq \frac{4}{3}$, then $d_{\Sigma}(p,q) \leq 4 < \rho(p)$ so let $\gamma$ be
	a shortest path in $\Sigma$ between $p$ and $q$. We set $L :=
		\|\gamma\|$. Let $\ell_1$,$\ell_2$ be the arclength-parametrized geodesics
	$o \to p$ and $o \to q$
	respectively. Clearly the concatenation $\tau := \ell_{2}^{-1} \circ \gamma
		\circ \ell_1$ is a loop of length $T := L + 2 r_*$. On this loop we define the
	piece-wise smooth unit-vector field $w(t)$ as $\ell'_1$ on $\ell_1$,
	$\ell'_2$ on $\ell_2$, and $\nu_\Sigma(\tau(t))$ on $\gamma$. Note that
	$w(0) \neq w(T)$.

	Let $P_{\tau, t\to 0}$ denote the parallel transport along $\tau$ from
	$\tau(t)$ to $\tau(0) = o$. We now consider the time-dependent vector $v(t)
		:= P_{\tau, t \to 0} w(t)$. On one hand we can write
	\begin{align*}
		v(T) - v(0) = (P_{\tau, T \to 0} - \Id) w(T) + (w(T) - w(0)).
	\end{align*}
	First let us bound the norm of $P_{\tau, T \to 0} - \Id$. Since $p$ and
	$q$ both lie on $\partial B$, we have $d_{N \times \RR}(p,q) \leq 2 r_*$.
	It follows that $L \leq \frac{4}{3} d_{N \times \RR}(p,q) \leq
		\frac{8}{3}r_*$. Since $r_* = 1$, this gives $r_* + L <
		4 \leq \inj(N)$, so the loop $\tau$ stays inside the injectivity ball
	centered at $o$. Let
	$\tilde{\tau} : [0,T] \times [0,1] \to N \times \RR$ be the homotopy between
	$\tau$ and the constant loop at $o$ given by
	\[
		\tilde{\tau}(t,s) = \exp_{o}(s \exp^{-1}(\tau(t))).
	\]
	We then have the area of this homotopy is given by
	\[
		\int_{[0,T] \times [0,1]} \|\partial_t \tilde{\tau} \wedge
		\partial_s \tilde{\tau} \| d t d s.
	\]
	Now along $\ell_1$ and $\ell_2$, we have that $\partial_t \tilde{\tau}$ and
	$\partial_s \tilde{\tau}$ are both radial. Hence for $t \in [0,r_*]$ and $t
		\in [r_*+L,2 r_*+L]$ the integrand vanishes. We thus have
	\[
		\int_{[0,T] \times [0,1]} \|\partial_t \tilde{\tau} \wedge
		\partial_s \tilde{\tau} \| d t d s
		=
		\int_{[r_*,r_*+L] \times [0,1]} \|\partial_t \tilde{\tau} \wedge
		\partial_s \tilde{\tau} \| d t d s
		\leq
		\|\gamma\| \max \| \partial_s \tilde{\tau} \|,
	\]
	where we used that for small enough curvature, the radial contraction is length-non increasing.
	The second factor here is bounded by the distance of $\gamma$ from $o$ hence
	is bounded by $r_* + L < 4$. Thus the area of the homotopy is at most $4 \|\gamma\|$.
	From here standard holonomy estimates, see \cite[\S6.2.1, p.~92]{alma993457963502466},
	show that the norm of $P_{\tau, T \to 0} - \Id$ is bounded above by $4 K
		\|\gamma\|$, which by the chord-arc estimate is bounded above by $6 K d_{N \times \RR}
		(p,q)$.

	Clearly we have that $r_* w(0) = \exp^{-1}_o (p)$ and $r_* w(T) =
		\exp^{-1}_o (q)$ so using Rauch comparison theorem we get $(1 - C K)
		d_{N \times \RR}(p,q) \leq r_* \|w(T) - w(0)\|$. We thus arrive at the
	estimate:
	\begin{align*}
		\|v(T) - v(0)\| \geq {}
		 & (1 - C K) r_*^{-1} d_{N \times \RR}(p,q) - 6 K d_{N
				\times \RR} (p,q) \geq \left(\frac{2}{3} - C K\right)
		d_{N \times \RR} (p,q)
		\\ \geq {}
		 & \left(\frac{2}{3} - C \alpha^2\right) d_{N \times \RR} (p,q).
	\end{align*}

	On the other hand we can also integrate to get
	\begin{align*}
		v(T) - v(0) = \int_{0}^T \partial_s v(s) d s = \int_{0}^T
		P_{\tau, t \to 0} (\nabla_{\tau'(t)} w(t)) d t.
	\end{align*}
	Since $\ell_1,\ell_2$ are geodesics, clearly $\nabla_{\tau'(t)}
		w(t)$ is zero on those segments. Along $\gamma$ we instead have
	$\nabla_{\tau'(t)} w(t) = A_{\Sigma} (\tau'(t))$. We thus have the estimate
	\[
		\|v(T) - v(0)\| \leq \int_{\gamma} \|A_{\Sigma} (\tau'(t))\| d t
		\leq L(\gamma) \kappa \leq \frac{4}{3} \kappa d_{N \times \RR}
		(p,q) \leq \frac{4}{3} \alpha d_{N \times \RR} (p,q).
	\]
	Combined with the previous estimate we get
	\[
		\frac{2}{3} - C \alpha^2 \leq \frac{4}{3} \alpha,
	\]
	for all sufficiently small $\alpha > 0$. This is clearly a contradiction and thus the theorem
	is proven.
\end{proof}

%% file: sections/examples.tex
\section{Examples}\label{sec:examples}
We show that graphicality and the
scale in \cref{lem:intro-chord-arc} are needed. All examples are in
Euclidean space.

\begin{example}[Hourglass]\label{ex-hourglass}
	For $R \geq 1$, let $\Sigma_R \subset \RR^n \times \RR$ be a smooth
	embedded hypersphere, flat on the two
	discs $\{|x| \leq R\} \times \{\pm 1\}$ and joining their boundaries by a
	rim of circular arcs of radius $R$ at distance of order $R$ from the axis.
	The profile curve and the rotation both contribute curvature at most $C/R$,
	so $|A_{\Sigma_R}| \leq C/R$. For the outer normal, the tangent ball at
	$(0,1)$ of radius $r > 1$ contains $(0,-1) \in \Sigma_R$, so $r_* \leq 1$
	there. As $R \to \infty$, $\Sigma_R$ converges to the two hyperplanes $\{z
	= \pm 1\}$; replacing $\pm 1$ by $\pm\varepsilon$ and letting $\varepsilon
	\to 0$ gives limits of multiplicity two.
\end{example}

\begin{example}[A tilted cylinder]\label{ex-cylinder}
	The bound in \cref{lem:intro-chord-arc} is sharp up to the constant $C$.
	Let $n \geq 2$, $0 < \varepsilon < 1/16$ and $L \geq 2$, and let
	\[
		u(x) := \frac{L}{2}\Bigl(\sqrt{1 + 8x_1 + x_2/(2\varepsilon)} - 1\Bigr)
	\]
	on $\Omega := \{x \in \RR^n : x_1 > -1/16,\ |x_2| < \varepsilon\}$. Take
	$x := 0$ and $y := e_1$, so that $r = 1$ and $B_\varepsilon(\gamma)
	\subseteq \Omega$. We show that $\delta \leq CL^{-2}$, while
	\[
		d_\Sigma(\hat x,\hat y) \geq (1 + c\delta)\,d(\hat x,\hat y).
	\]

	The graph $\Sigma$ of $u$ is part of the parabolic cylinder $\{F = 0\}$,
	where $F := x_2 + 16\varepsilon x_1 - 8\varepsilon L^{-2}t(t+L)$. Since
	$|\nabla F| \geq 1$ and $|\Hess F| = 16\varepsilon L^{-2}$, we have
	$\kappa \leq C\varepsilon L^{-2}$, so $\delta = \kappa/\varepsilon \leq
	CL^{-2}$. Let $g(s) := \sqrt{s^2 + L^2/4} + \sqrt{(1-s)^2 + L^2/4}$ be the
	length of the broken line from $\hat x = (0,0)$ to $\hat y = (e_1, L)$
	through $(se_1, L/2)$; then $d(\hat x,\hat y) = g(1/2)$. A path in $\Sigma$
	from $\hat x$ to $\hat y$ meets $\{t = L/2\}$ at a point $(z, L/2)$, where
	$u(z) = L/2$ forces $z_1 = 3/8 - z_2/(16\varepsilon) < 7/16$. Its
	projection to the $(x_1,t)$-plane is no shorter, so the path has length at
	least $g(z_1) \geq g(7/16)$, since $g$ is convex and decreasing on
	$[0,1/2]$. Now $g'(1/2) = 0$ and $g'' \geq c/L$ on $[0,1]$, so $g(7/16)
	\geq g(1/2) + c/L$; dividing by $g(1/2) \leq CL$ gives $d_\Sigma(\hat
	x,\hat y)/d(\hat x,\hat y) - 1 \geq cL^{-2} \geq c\delta$.
\end{example}

%% file: sections/applications.tex
\section{Applications}\label{sec:apps}
We present two applications of these results. First we prove
an in-radius analog of the Heinz curvature estimate for
minimal graphs.
\subsection*{A Heinz-type estimate}
\begin{corollary}\label{cor-heinz-inscribed-curvature}
	For $n \leq 7$, if $\Sigma \seq \RR^{n+1}$ is a minimal graph over a
	Euclidean ball $D_{r_0}\seq \RR^n$ and $0 < r < r_0$ then,
	\[
		\sup_{D_{r_0 - r}} r_*^{-1} \leq C(n) r^{-1}.
	\]
\end{corollary}
\begin{proof}
	By the standard Heinz curvature estimate for minimal graphs in dimensions
	$n\leq7$, see \cite[Theorem~2.24, p.~88]{coldingCourseMinimalSurfaces2011},
	we know that
	\[
		\sup_{D_{r_0 - r/2}} \|h_\Sigma\| \leq 2 C r^{-1}.
	\]
	For $p$ over $D_{r_0-r}$, any tangent-ball obstruction coming from outside
	$D_{r_0-r/2}$ has radius at least $r/4$.
	Such an obstruction contributes at most $4 r^{-1}$ to $r_*^{-1}$.
	Otherwise we apply \cref{thrm-inscribed-curvature-estim} to the restricted
	graph over $D_{r_0-r/2}$ and get
	\[
		r_*(p)^{-1} \leq C \kappa + 4 \rho(p)^{-1}
		\leq C r^{-1} + 4 (r/2)^{-1} \leq C r^{-1}.
	\]
\end{proof}

\subsection*{Area estimates}
Next we prove uniform area estimates depending only on extrinsic
curvature for graphs over any Riemannian manifold. One estimate
follows from \cref{lem-chord-arc} and the other from
\cref{thrm-inscribed-curvature-estim}.

\begin{corollary}\label{cor-local-area-estimate}
	With $\kappa,K$ as in \eqref{eq-kappa} and
	\eqref{eq-ambient-curvature}, there exist dimensional constants $c, C >
		0$ such that for any $p \in \Sigma$, if
	\[
		r < \min(c \kappa^{-1}, c K^{-1/2}, \inj(N), d(p,\partial \Omega)/2),
	\]
	then
	\[
		|B_{r}^{N \times \RR}(p) \cap \Sigma|
		\leq
		\omega_n r^n (1 + C \kappa r + C K r^2).
	\]
\end{corollary}
\begin{proof}
	Let $x := \pi(p)$ and let $q \in B_r^{N \times \RR}(p) \cap \Sigma$ lie
	over $y \neq x$, so $d(x,y) < r$. We take a minimizing geodesic $\gamma$
	of $N$ from $x$ to $y$ and $\varepsilon := d(x,y)/3$. Then
	$B_\varepsilon(\gamma) \seq B(x,2r) \seq \Omega$, since $2r < d(x,\partial
	\Omega)$, and $\delta = 3\kappa d(x,y) + K d(x,y)^2 \leq 3\kappa r + Kr^2 < c$,
	which is below $\delta_0$ for $c$ small. So \cref{lem-chord-arc} gives
	$d_\Sigma(p,q) \leq (1+C\delta) d_{N \times \RR}(p,q)$, and hence
	\[
		B_r^{N \times \RR} (p) \cap \Sigma \seq B_s^\Sigma (p), \qquad
		s := r (1 + C(\kappa r + K r^2)) \leq 2r.
	\]
	By the Gauss equation, $\Ric_\Sigma \geq -C(K+\kappa^2)g_\Sigma$. Since
	$\pi$ does not increase distances, the closed intrinsic ball of radius $2r$
	about $p$ lies over $\ov{B(x,2r)} \seq \Omega$, so it is compact. We may
	therefore apply Bishop--Gromov comparison, which gives
	\[
		|B_s^\Sigma(p)| \leq \omega_n s^n \bigl(1 + C(K+\kappa^2)s^2\bigr).
	\]
	Expanding $s^n$ and using $s \leq 2r$ and $\kappa^2r^2\leq c\kappa r$
	gives the result.
\end{proof}
\begin{corollary}\label{cor-uniform-area-estimate}
	With $\kappa,K$ as in \eqref{eq-kappa} and
	\eqref{eq-ambient-curvature}, for any $R > r > 0$, there exists a constant
	$C = C(\kappa, K, r, R, \inj (N))$ such that for any $q \in N \times \RR$
	with $R \leq d(q,\partial \Omega)$ we have
	\[
		|B_r^{N \times \RR} (q) \cap \Sigma| \leq C.
	\]
\end{corollary}
\begin{proof}
	Let $\Sigma_r := B_r^{N\times\RR}(q)\cap\Sigma$. Points of $\Sigma_r$ lie
	over $B(\pi(q),r)$, so their distance to $\partial\Omega$ is more than
	$R - r$. By \cref{thm:main} for both orientations, both in-radii at points
	of $\Sigma_r$ are at least
	\[
		\rho_0 := \bigl(C(\kappa + K^{1/2}) + 4\min(R-r,\inj(N))^{-1}\bigr)^{-1},
	\]
	hence $(p,t) \mapsto
	\exp_p(t\nu(p))$ is injective on $\Sigma_r \times (-\rho_0,\rho_0)$. Its
	Jacobian $J$ satisfies $|\partial_t \log J| = |H|$, where $H$ is the mean
	curvature of the parallel hypersurface at distance $t$, and since
	$\|h_\Sigma\| \leq \kappa$ and $\|\Rm_N\| \leq K$, Riccati comparison
	gives $|H| \leq C(\kappa + K^{1/2})$ for
	$|t| \leq c(\kappa + K^{1/2})^{-1}$. So with $\rho_1 := \min(\rho_0,
	c(\kappa + K^{1/2})^{-1})$ we get $J \geq e^{-C}$ for $|t| < \rho_1$, and
	since the image of $\Sigma_r \times (-\rho_1,\rho_1)$ lies in
	$B_{r+\rho_1}^{N\times\RR}(q)$,
	\[
		c\rho_1 |\Sigma_r| \leq |B_{r+\rho_1}^{N\times\RR}(q)|.
	\]
	Volume comparison bounds the right hand side in terms of $K$, $r$ and
	$\rho_1$, and $\rho_1$ depends only on $\kappa,K,r,R$ and $\inj(N)$.
\end{proof}

%% file: sections/app_differential_geometry.tex
\appendix
\section{Comparison Principles and Metric Estimates}\label{sec:app-metric-comparison}

Let $M_1,M_2$ be two hypersurfaces in a Riemannian manifold $N$
	which intersect at the point $p$ with $\nu_{M_1} (p) = \nu_{M_2} (p)$.
	We say that $M_1$ is locally interior to $M_2$ if there is some
	neighborhood $U$ of $p$ such that $U \cap (M_2 \setminus M_1)$ lies
	entirely in the connected component corresponding to $\nu_{M_1} (p)$.
\begin{proposition}\label{prop-onesided}
	If $M_1$ is locally interior to $M_2$ then the
	second fundamental forms satisfy
	\begin{align}
		h_{M_2} (p) \leq h_{M_1} (p),
	\end{align}
	in the sense of quadratic forms on the tangent space $T_p M_1 = T_p
		M_2$. Similarly if $h_{M_2} (p) < h_{M_1} (p)$ then $M_1$ is
	locally interior to $M_2$.
\end{proposition}
\begin{proof}
	Choose a local defining function $f$ for $M_1$ normalized by
	$\ov{\nabla}_{\nu_{M_1}(p)}f=1$. If $M_1$ is locally interior to
	$M_2$, then $f|_{M_2}\geq0$ near $p$. Restricting $f$ to any geodesic in
	$M_2$ through $p$ and applying the second-derivative test gives
	$h_{M_2}(v,v)\leq h_{M_1}(v,v)$ for every unit tangent vector $v$.
	Conversely, if $h_{M_2}<h_{M_1}$ as quadratic forms, the same Taylor
	expansion gives $f|_{M_2}>0$ away from $p$ in every tangent direction. By
	compactness of the unit sphere this holds in a full neighborhood, so
	$M_1$ is locally interior to $M_2$.
\end{proof}

\begin{proposition}[Almost linear tube coordinates]\label{prop-connection-expansion}
	Fix $\ell > \varepsilon > 0$, and let $\gamma : [-\varepsilon,\ell+\varepsilon] \to M$ be a geodesic in an
	$n$-dimensional complete Riemannian manifold with distance at least $\varepsilon$ from $\partial M$. Set $U_\varepsilon := [-\varepsilon, \ell + \varepsilon] \times B_{\varepsilon}^{n-1}$.
	There exist dimensional constants $c,C > 0$ independent of $\ell$ such that if
	\[
		\varepsilon^2 K \leq c,\qquad
		K := \|\Rm_{g}\|_{C^0(M,g)},
	\]
	then there is a map
	$\Phi : U_{c\varepsilon} \to M$ with
	\[
		\Phi(0) = \gamma(0),\qquad \Phi(\ell e_1) = \gamma(\ell),
	\]
	with $\Phi(U_{c\varepsilon}) \seq B_{\varepsilon} (\gamma)$. Moreover, the pullback $\tilde{g} = \Phi^* g$ satisfies
	\[
		\|\tilde{g} - \delta\| \leq C\varepsilon^2 K,\qquad
		\|\partial \tilde{g}\| \leq C\varepsilon K.
	\]
	In particular, if $\ov{\nabla}$ is the Euclidean connection in these
	coordinates, then
	\[
		\|\nabla^{\tilde{g}} - \ov{\nabla}\|_{C^0(U_{c\varepsilon},\tilde{g})}
		\leq C\varepsilon K.
	\]
\end{proposition}
\begin{proof}
	After rescaling the metric by $\varepsilon^{-2}$, it is enough to prove
	the result for $\varepsilon = 1$. Set $\theta = K$ in this scale. At the
	end of the proof we scale back, replacing $\theta$ by
	$\varepsilon^2 K$.

	Since $\|\Rm_g\| \leq \theta \leq c$, the Fermi coordinate map $F$ of $\gamma$ is an immersion on $U_1$. By first pulling back by $F$, we may replace $(M,g)$ by $(U_1, F^*g)$, and thus get by Jacobi comparison that
	\[
		\|g-\delta\|\leq C \theta.
	\]
	Fermi coordinates do not control $\partial g$ without a bound on $\nabla \Rm$, so we pass to harmonic coordinates.
	Choose a dimensional constant $a > 0$ and points $p_j = \gamma(s_j)$ with fixed spacing, so that the Euclidean balls $B_\delta(p_j,2a)$ cover a neighborhood of $U_a$ with uniformly bounded overlap and $B_j(3a) \seq U_1$, where $B_j(r) = B_g(p_j,r)$.
	On each $B_j (3 a)$ choose harmonic coordinates $u_j : B_j (3 a) \to \RR^n$ with
	almost linear Dirichlet boundary data, and postcompose by an affine map so that
	\[
		u_j(p_j) = p_j = s_j e_1 ,\qquad Du_j|_{p_j} = \Id.
	\]
	The estimates in \cite[Equation 2.8.10]{jostHarmonicMappings1984} compare gradients of harmonic maps with radially transported vectors. Since Fermi coordinates are constructed from parallel transports, these estimates can be used to give $D u_j = \Id + O(\theta)$ on fixed interior subballs of $B_j(3a)$. Because each ball is anchored at $p_j$, in the sense $p_j$ is fixed, integrating from $p_j$ gives $u_j = \Id + O(\theta)$ there as well, importantly this is independent of $\ell$.
	For small enough $\theta$, the coordinate images therefore contain $B_\delta(p_j,2a)$. Define the pushforward metrics $h_j = (u_j)_* g$. The quantitative estimates in \cite[Lemma 2.8.2]{jostHarmonicMappings1984} give on these coordinate balls
	\[
		\|h_j-\delta\|+\|\partial h_j\|\leq C \theta.
	\]
	In particular, the Christoffel symbols of each $h_j$ satisfy
	\[
		\|\Gamma(h_j)\|\leq C \theta.
	\]
	On an overlap of these restricted charts, put $\tau = u_k\circ u_j^{-1}$. The exact coordinate-change formula for Christoffel symbols gives at any point $q$ in its coordinate domain,
	\[
		\partial_a\partial_b\tau^\mu|_{q}
		=
		\Gamma^c_{ab} (h_j)|_{q} \partial_c\tau^\mu|_q
		-
		\Gamma^\mu_{\alpha\beta}(h_k)|_{\tau (q)}
		\partial_a\tau^\alpha|_q
		\partial_b\tau^\beta|_q.
	\]
	The estimates on $D\tau$ and the Christoffel symbols give
	$\|\tau-\Id\|_{C^2} \leq C\theta$ on each coordinate overlap.

	Choose a smooth partition of unity $\{\eta_j\}$ on a neighborhood of $U_a$,
	subordinate to the balls $B_\delta(p_j,2a)$, with uniformly bounded first
	and second derivatives. Extend $\eta_j\circ u_j$ by zero outside its
	chart, and define on $U_a$
	\[
		\chi_j=\frac{\eta_j\circ u_j}{\sum_k\eta_k\circ u_k},
		\qquad \Psi=\sum_j\chi_j u_j.
	\]
	The denominator is $1+O(\theta)$, and the transition estimates give
	uniform $C^2$ bounds for the $\chi_j$ in every harmonic chart. Writing
	$\tau_{jk}=u_k\circ u_j^{-1}$, we have in the $u_j$ chart
	\[
		Q_j:=\Psi\circ u_j^{-1}
		=\Id+\sum_k(\chi_k\circ u_j^{-1})(\tau_{jk}-\Id).
	\]
	Thus, independently of $\ell$,
	\[
		\|Q_j-\Id\|_{C^2}\leq C\theta,
		\qquad \|\Psi-\Id\|_{C^1(U_a)}\leq C\theta.
	\]
	Since $\ell>1$, the endpoint errors allow an affine map $A$ with
	\[
		A(\Psi(0))=0,\qquad A(\Psi(\ell e_1))=\ell e_1,
		\qquad \|DA-\Id\|\leq C\theta/\ell.
	\]
	Then $\|A\circ\Psi-\Id\|_{C^1(U_a)}\leq C\theta$, so for
	$c$ small enough this map is bi-Lipschitz and its image contains
	$U_c$. Set
	$\Phi=(A\circ\Psi)^{-1}$ on $U_c$. In each chart,
	$\Phi^*g=((A\circ Q_j)^{-1})^*h_j$, so the preceding $C^2$ bounds
	and the harmonic estimates on $h_j$ give
	\[
		\|\Phi^* g - \delta\|+\|\partial(\Phi^* g)\|
		\leq C \theta.
	\]
	Composing with the initial Fermi map returns to the original manifold.
	Scaling back gives
	\[
		\|\Phi^* g-\delta\| \leq C\varepsilon^2K,\qquad
		\|\partial(\Phi^* g)\| \leq C \varepsilon K,
	\]
	and the connection estimate follows from the formula for the Christoffel
	symbols.
\end{proof}

\begin{proposition}\label{prop-h-bounds}
	Let $\iota : \Sigma \to M$ be an immersion of a hypersurface. Let
	$g_1$ and $g_2$ be two metrics on $M$, with connections $\nabla^1$
	and $\nabla^2$, such that $\|g_1 - g_2\|_{g_1} < \frac{1}{2}$. Denote by $h_1$ and
	$h_2$ the second fundamental forms of $\Sigma$ with respect to
	$g_1$ and $g_2$, respectively.
	Then we have the estimate
	\[
		\|h_2\|_{g_2} \leq C(n)(\|h_1\|_{g_1}
		+ \|\nabla^1 - \nabla^2\|_{g_1}),
	\]
	at any point of $\Sigma$.
\end{proposition}
\begin{proof}
	Fix a $g_1$-orthonormal frame adapted to $\Sigma$, and let $\omega$ be
	the corresponding conormal. For both metrics, $h_i(X,Y)$ can be written
	as $-\|\omega\|_{g_i}^{-1}\omega(\nabla^i_XY)$. Since the norms of $g_1$ and $g_2$ are uniformly equivalent up to a dimensional factor, we can ignore the norm factors and unit-vector normalizations. Replacing
	$\nabla^2_X X$ by $\nabla^1_X X$ in the formula above contributes at most
	$\|\nabla^1-\nabla^2\|_{g_1}$ up to a dimensional constant. Taking the supremum of
	$|h_2(X,X)|$ over $g_2$ unit vectors $X$ then gives the proposition.
\end{proof}

%% file: sections/app_analytic_estimates.tex
\section{Auxiliary Lemmas}\label{sec:app-auxiliary}

\begin{proposition}\label{prop-qinverse}
	Let $B_r$ be a complete Riemannian ball with boundary of radius $r$, and
	let $f : B_r \to N$ be a $C^1$ map to a connected Riemannian manifold
	such that for some $a > 0$, $B_{a} \seq Df_p (B_1)$ for
	all $p \in B_r$. Then the image $f(B_r)$ contains $B_{a r} (f(0))$.
\end{proposition}
\begin{proof}
	Given $z \in B_{ar}(f(0))$, choose $\delta > 0$ with $d(f(0),z) <
	a(r-\delta)$, and join $f(0)$ to $z$ by a path on $[0,1]$ of speed at most
	$a(r-\delta)$. The set of times for which this path has a
	horizontal lift starting at $0$ is nonempty and open by local horizontal
	lifting. The assumption $B_a\subset Df_p(B_1)$ gives a horizontal preimage of
	each velocity vector with norm at most $r-\delta$, so every partial lift has
	length at most $r-\delta$ and remains in the compact subball
	$\overline{B_{r-\delta}}$. Thus the set of times is also closed, and thus we get the statement by continuity.
\end{proof}

%% file: sections/app_geometric_identities.tex
\begin{proposition}\label{prop-simple-geometric}
	Let $v$ and $w$ be unit vectors in $\RR^{n}$. Then the following holds.
	\begin{enumerate}[label=(\roman*)]
		\item $\ip{v,w} = 1 - \|v - w\|^2/2$.
		\item $\|\pi_{w^\perp} v\| \leq \|v - w\|$.
		\item If $\|v - w\| \leq 1$ then $\|\pi_{w^\perp} v\| \geq
			      \frac{\sqrt{3}}{2}\|v - w\|$.
		\item If $\pi_{w^\perp} v \neq 0$ then $\|v - \frac{\pi_{w^\perp} v}{\|\pi_{w^\perp} v\|}\| \leq
			      \sqrt{2} |\ip{v,w}|$.
	\end{enumerate}
\end{proposition}
\begin{proof}
	Decompose $v$ into its components parallel and orthogonal to $w$ and use
	$\|v-w\|^2=2(1-\ip{v,w})$ which gives $(i)$. Projections are one Lipschitz which gives $(ii)$.  Then
	$\|\pi_{w^\perp}v\|^2=1-\ip{v,w}^2
		=\|v-w\|^2(1-\|v-w\|^2/4)$ gives $(iii)$ immediately. Finally,
	\[
		\left\|v-\frac{\pi_{w^\perp}v}{\|\pi_{w^\perp}v\|}\right\|^2
		=2-2\|\pi_{w^\perp}v\|
		\leq 2\ip{v,w}^2,
	\]
	which gives $(iv)$.
\end{proof}